\documentclass[a4paper,reqno,10pt,oneside]{amsart}
\usepackage{amsmath,geometry}
\usepackage{graphicx}
\usepackage{mathrsfs}
\usepackage{amssymb,amsmath, amsfonts, amsthm}
\usepackage{latexsym}
\usepackage{color} %{\color{red}..............OK}
\usepackage[dvips]{epsfig}
\usepackage{graphicx}

\allowdisplaybreaks[1]

\numberwithin{equation}{section}

\renewcommand{\(}{\left(}
\renewcommand{\)}{\right)}
\renewcommand{\[}{\left[}
\renewcommand{\]}{\right]}

\newtheorem{theorem}{Theorem}[section]
\newtheorem{proposition}[theorem]{Proposition}

\newtheorem{lemma}[theorem]{Lemma}
\newtheorem{remark}[theorem]{Remark}

\newcommand{\F}{{\cal F}}

\newcommand{\beq }{\begin{equation}}
\newcommand{\eeq }{\end{equation}}

\def\F{{\mathcal F}}

\newcommand{\beqs}{\begin{equation*}}
\newcommand{\eeqs}{\end{equation*}}
\newcommand{\beqn}{\begin{eqnarray}}
\newcommand{\eeqn}{\end{eqnarray}}
\newcommand{\beqns}{\begin{eqnarray*}}
\newcommand{\eeqns}{\end{eqnarray*}}
\newcommand{\bdoc}{\begin{document}}
\newcommand{\edoc}{\end{document}}
\newcommand{\be}{\begin{enumerate}}
\newcommand{\ee}{\end{enumerate}}
\newcommand{\bdescr}{\begin{description}}
\newcommand{\edescr}{\end{description}}
\newcommand{\ba}{\begin{array}}
\newcommand{\ea}{\end{array}}
\newcommand{\intR}{\int_{\mathbb R^N}}
\newcommand{\R}{\mathbb R}
\newcommand{\Mm}{\mathcal{M}_{\lambda,\Lambda}^-}
\newcommand{\Mp}{\mathcal{M}_{\lambda,\Lambda}^+}
\newcommand{\Mpm}{\mathcal{M}_{\lambda,\Lambda}^\pm}
\newcommand{\Nm}{\tilde N_-}
\newcommand{\Np}{\tilde N_+}
\newcommand{\Nmp}{\tilde N_\mp}
\newcommand{\parallelsum}{\mathbin{\!/\mkern-5mu/\!}}
\newcommand{\e}{\varepsilon}
 \renewcommand{\(}{\left(}
\renewcommand{\)}{\right)}
\renewcommand{\[}{\left[}
\renewcommand{\]}{\right]}
\newenvironment{Proof}{\noindent{\bf Proof}}{\hfill$\Box$\\[2mm]}

\begin{document}
\title[Liouville-type results]{Liouville-type results in exterior domains for radial solutions  of fully nonlinear equations}

\author{G\MakeLowercase{iulio} Galise, A\MakeLowercase{lessandro} Iacopetti, F\MakeLowercase{abiana} Leoni}

\thanks{\emph{Acknowledgements.} This research is partially supported by INDAM-GNAMPA}
\address{Dipartimento di Matematica, Sapienza Universit\`a di Roma,  P.le Aldo Moro 2, 00185 Roma, Italy}
\email{galise@mat.uniroma1.it (G. Galise)}
\email{iacopetti@mat.uniroma1.it (A. Iacopetti)}
\email{leoni@mat.uniroma1.it (F. Leoni)}

\begin{abstract}
We give necessary and sufficient conditions for the existence  of positive radial solutions for a class of fully nonlinear uniformly elliptic equations posed in the complement of a ball in $\R^N$, and equipped with homogeneous Dirichlet boundary conditions.
\end{abstract}

\maketitle

\section{Introduction}
Let $B$ be any ball in $\R^N$  and let $0<\lambda\leq \Lambda$. The aim of this paper is to detect the optimal conditions  on  the exponent $p>1$ for the existence of positive radial solutions of the fully nonlinear exterior Dirichlet problem
\begin{equation}\label{eq:probMpm}
\begin{cases}
-\mathcal{F}(D^2 u)=u^{p} & \hbox{in} \ \R^N\setminus \overline B,\\
\qquad  \ \ \ \ \ u=0 & \hbox{on} \ \partial B,
\end{cases}
\end{equation}
where $\mathcal{F}$ is either one of the Pucci's extremal operators   $\Mpm$, defined respectively as
%$$
%\begin{array}{c}
%{\mathcal M}^-_{\lambda,\Lambda}(X)= \displaystyle{\inf_{ \lambda I\leq A\leq \Lambda I}{\rm tr} (AX)=\lambda\, \sum_{\mu_i>0}\mu_i+\Lambda\, \sum_{\mu_i<0}\mu_i} \\[2ex]
%{\mathcal M}^+_{\lambda,\Lambda}(X)= \displaystyle{\sup_{ \lambda I\leq A\leq \Lambda I} {\rm tr} (AX)=\Lambda\, \sum_{\mu_i>0}\mu_i+\lambda\, \sum_{\mu_i<0}\mu_i}
%\end{array}
%$$
\begin{equation*}
\begin{split}
{\mathcal M}^-_{\lambda,\Lambda}(X)&:=\inf_{ \lambda I\leq A\leq \Lambda I}{\rm tr} (AX)=\lambda\, \sum_{\mu_i\geq 0}\mu_i+\Lambda\, \sum_{\mu_i<0}\mu_i \\
{\mathcal M}^+_{\lambda,\Lambda}(X)&:=\sup_{ \lambda I\leq A\leq \Lambda I} {\rm tr} (AX)=\Lambda\, \sum_{\mu_i\geq 0}\mu_i+\lambda\, \sum_{\mu_i<0}\mu_i
\end{split}
\end{equation*}
 $\mu_1,\ldots,\mu_N$ being the eigenvalues of any squared symmetric matrix $X$.

Pucci's extremal operators are the prototype of fully nonlinear uniformly elliptic operators.  Acting  as barriers in the whole class of  operators with fixed ellipticity constants $\lambda\leq \Lambda$,  they
play a crucial role in the  regularity theory for fully nonlinear elliptic equations, see \cite{CC}. Moreover, as $\sup$/$\inf$ of linear operators, they frequently arise in  the equations satisfied by the value function associated with stochastic optimal control problems, see e.g. \cite{FS, L}, with special application to mathematical finance problems.

We recall that if $\Lambda=\lambda$, both Pucci's operators reduce, up to a multiplicative factor, to the Laplace operator. Thus, they may be considered also as perturbations of the standard Laplacian, and  the well-known Lane-Emden-Fowler equation $-\Delta u=u^p$ is included as a very special case of the problems we are considering.\\
 In the semilinear case,  the Lane-Emden-Fowler equation has been largely studied. The well known existence results, exhibiting the critical Sobolev exponent $p^*=\frac{N+2}{N-2}$ as threshold for the existence of entire solutions or solutions in bounded domains,   are intimately related to the  (lack of) compactness properties of  the Sobolev embeddings. Moreover,  the entire solutions existing in  the critical case $p=p^*$ realize the best constant in the Sobolev inequality, see \cite{T}.  In other words, the critical  nature of the exponent $p^*$ may be largely  interpreted  in view of  structural  properties of the functional setting behind the equation. 

As soon as $\Lambda>\lambda$, Pucci's operators loose the linear and  variational structures. Nevertheless, Lane-Emden-Fowler type equations as \eqref{eq:probMpm} have been studied and, at least in the radial setting, some critical exponents $p^*_\pm$ acting as thresholds for the existence of entire radial solutions or solutions in balls have been proved to exist, see \cite {FQ}. In the fully nonlinear radial setting, the exponents $p^*_\pm$ play the same role as the critical Sobolev exponent $p^*=\frac{N+2}{N-2}$ for the Laplacian. Though their appearance is motivated exclusively as threshold for  the existence of entire radial solutions or solutions in balls, the recent results of \cite{BGLP}, where some  weighted energies associated with radial solutions of \eqref{eq:probMpm} are introduced and proved to be asymptotically preserved by almost critical solutions, suggest that the critical exponents   reflect some intrinsic properties of the operators,  maybe  beyond the radial setting.\\
The Liouville-type results obtained in the present paper may be regarded as a further justification of the critical character of the exponents $p^*_\pm$, since they are proved to act as thresholds also for the existence of radial solutions to exterior Dirichlet problems. We observe that this result, while well expected for  semilinear equations  due to the duality between Dirichlet problems in balls and Dirichlet problems in exterior domains (via the Kelvin transform), in the fully non linear case it requires a direct proof.
Moreover,  the existence  of positive solutions for exterior Dirichlet problems is also closely related to the existence or non-existence of sign-changing radial solutions in balls or in the whole space (see also Remark 5.2 in \cite{GLP}). In particular, in the recent paper \cite{GILP},    the asymptotic  behavior of sign changing solutions in balls  is studied, and the obtained results strongly rely on the present  theorems for exterior Dirichlet problems.

In order to describe the results, let us introduce  the dimension like parameters
$$
\begin{array}{c}
\tilde{N}_-\, :\, = \frac{\Lambda}{\lambda} (N-1)+1\quad \hbox{ for } \Mm\, ,\\[2ex]
\tilde{N}_+\, :\, = \frac{\lambda}{\Lambda} (N-1)+1\quad \hbox{ for } \Mp\, ,
\end{array}
$$
which have been proved in some previously studied cases to play a key role in existence results.   
Let us emphasize  that one has always $\tilde{N}_+\leq N\leq \tilde{N}_-$, and the equalities hold true if and only if $\Lambda=\lambda$.

The case of \emph{entire supersolutions} has been   considered in \cite{CL} ,  where it has been proved that
$$
\exists\, u >0\, ,\ -\Mpm (D^2u)\geq u^p\quad \hbox{ in } \R^N \Longleftrightarrow p >\frac{\tilde{N}_\pm}{\tilde{N}_\pm-2}\, ,
$$
meaning that, in particular, positive supersolutions never exist if $\tilde{N}_\pm\leq 2$. In the sequel, we will always assume that $\tilde{N}_\pm>2$. 

The same threshold has been proved in \cite{AS} to be optimal for the existence of \emph{solutions in any exterior domain}, that is
$$
\exists\, u >0\, ,\ -\Mpm (D^2u)= u^p\quad \hbox{ in } \R^N\setminus K \Longleftrightarrow p >\frac{\tilde{N}_\pm}{\tilde{N}_\pm-2}\, ,
$$
where $K\subset \R^N$ is any nonempty compact set. Let us also mention the results of \cite{AGQ}, where more general nonlinearities $f(u)$ replacing $u^p$ are considered. In the above results, supersolutions are meant in the viscosity sense and no symmetry property on  $u$ is required.

In the present paper we are concerned with \emph{solutions}  of the equation satisfying further \emph{homogeneous Dirichlet boundary conditions}. By elliptic regularity theory, it is not restrictive to consider classical solutions of problem \eqref{eq:probMpm}, that are $C^2$ functions satisfying pointwise the equation as well as  the boundary  condition. In its full generality, that is without assuming radial symmetry of the solutions, the problem is completely open, and also in the semilinear case (i.e. when $\Lambda =\lambda$) few results are known for solutions of exterior Dirichlet problems, see e.g. \cite{DdPMW} where  solutions are constructed as perturbations of   radial  solutions. Our results, limited to  radially symmetric  solutions, may hopefully contribute to tackle the general fully nonlinear problem.

In the radial setting, the existence of \emph{entire positive  solutions} has been studied in \cite{FQ}, where  it has been proved that there exist two critical exponents $p^*_+$ and $p^*_-$ associated with $\Mp$ and $\Mm$ respectively, such that
$$
\exists\, u  \hbox{ radial}, u>0\, ,\ -\Mpm (D^2u)= u^p\quad \hbox{ in } \R^N  \Longleftrightarrow p \geq p^*_\pm\, .
$$
Unfortunately, the  dependence of the radial critical exponents on the effective dimensions is not explicitly known. The radial critical exponents  are proved in \cite{FQ}  to  satisfy, when $\lambda<\Lambda$, the strict inequalities
\begin{equation}\label{criticalexp}
\frac{\tilde{N}_-+2}{\tilde{N}_--2}< p^*_- < \frac{N+2}{N-2}< p^*_+< \frac{\tilde{N}_++2}{\tilde{N}_+-2}\,.
\end{equation}
Note that inequalities in \eqref{criticalexp} become equalities when $\Lambda=\lambda$.  Let us emphasize that inequalities \eqref{criticalexp} say that,  with respect to the intrinsic dimensions, the exponent  $p^*_+$  is subcritical, whereas $p^*_-$ is  supercritical.  Moreover, \eqref{criticalexp} show that the critical Sobolev exponent $p^*=\frac{N+2}{N-2}$ is not preserved even for small perturbations of the  operator, since the critical exponents are  instantaneously different  as soon as $\Lambda>\lambda$. For an integral characterization of $p^*_\pm$, as well as for sharp estimates on entire critical  solutions, we refer to  \cite{BGLP}.

Clearly, the analysis on the existence of entire positive radial solutions yields, as a by product, the dual result on the existence of positive solutions of Dirichlet problems in balls, namely
$$
\exists\, u>0\, ,\ -\Mpm (D^2u)= u^p\quad \hbox{ in } B, \ u=0 \ \hbox{ on } \partial B  \Longleftrightarrow p < p^*_\pm\, .
$$
Note that, in this case, the radial symmetry of the solutions is not a restriction, since, by \cite{DLS}, any positive solution in the ball is radial. 
The critical exponents $p^*_\pm$, therefore,  give the optimal thresholds for the existence of positive solutions in balls and, as proved in \cite{EFQ2}, also in domains sufficiently close to balls.

On the other hand, as recently proved in \cite{GLP}, for annular domains (radial) solutions exist for any $p>1$.  More precisely, in \cite{GLP} it has been proved that solutions of the initial value problems for the ODEs associated with the equations $-\Mpm(D^2u)=u^p$ give radial solutions in   annular domains provided that they have sufficiently large initial slope.
 The results of the present paper, in a sense, complement the results of \cite{GLP}, since here we  prove, in particular,  that a sufficiently large  initial velocity  is needed for having radial solutions in  annuli if and only if $p>p^*_\pm$.

The results of the subsequent  sections are summarized in the following theorem.

\begin{theorem}\label{main} There exist positive radial solutions of problem \eqref{eq:probMpm} if and only if $p>p^*_\pm$. Moreover, for any $p>p^*_\pm$, problem \eqref{eq:probMpm} has a unique positive radial solution $u^*$ satisfying 
\begin{equation}\label{fast}
\lim_{r\to +\infty} r^{\tilde{N}_\pm-2}u^*(r)=C>0\, ,
\end{equation}
and infinitely many positive radial solutions $u$ satisfying 
\begin{equation}\label{slow}
\lim_{r\to +\infty} r^{\tilde{N}_\pm-2}u(r)=+\infty\, .
\end{equation}
\end{theorem} 
Borrowing the terminology currently used, the solution $u^*$ satisfying \eqref{fast} will be referred to as the \emph{fast decaying} solution. As far as solutions $u$ satisfying \eqref{slow} are concerned, they will be proved to satisfy either
$$
\lim_{r\to +\infty} r^{\frac{2}{p-1}}u(r)=c>0\, ,
$$
in which case they will be called \emph{slow decaying} solutions, or
$$
0<\liminf_{r\to +\infty} r^{\frac{2}{p-1}}u(r)<\limsup_{r\to +\infty} r^{\frac{2}{p-1}}u(r)<+\infty\, ,
$$
in which case they will be named \emph{pseudo-slow decaying} solutions (see \cite{FQ}).
\medskip

In the semilinear case, a proof of Theorem \ref{main} can be found in \cite{JPY}. For the fully nonlinear case, it will be a straightforward consequence of the results proved in the next sections, where we perform  a careful analysis of the initial value problem for the ODE associated to radial solutions of problem \eqref{eq:probMpm}. Recalling that the eigenvalues of the Hessian matrix $D^2u$ of a smooth radial function $u=u(r)$ are nothing but $u''(r)$ and $\frac{u'(r)}{r}$ (with multiplicity at least $N-1$), it is not difficult to write  the ODE satisfied by a radial solution of problem  \eqref{eq:probMpm}. Nevertheless, since the coefficients of the operators $\Mpm$ depend   on the sign of the eigenvalues of the Hessian matrix, we will obtain an ODE with discontinuous coefficients having jumps at the points where the solution $u$ changes its monotonicity and/or concavity. This is a feature of the fully nonlinear problem which makes  techniques previously developed  for the semilinear case  not directly applicable. In particular,  the Kelvin transform which reduces a supercritical exterior Dirichlet problem to a subcritical  Dirichlet problem in the punctured ball cannot be used.

We will essentially make use of  the results of \cite{FQ} for  entire solutions, in particular of the fact that the critical exponents are  the only exponents for which the entire solutions are fast decaying.  Moreover, as in \cite{FQ}, we will take advantage of the Emden-Fowler trasformed
$$
x (t)=r^{\frac{2}{p-1}} u(r), \ \ \ r=e^t,
$$
which produces a new variable $x(t)$ satisfying an autonomous equation. Despite the fact that also the coefficients of the equation satisfied by $x$ will have jumps at the points corresponding to the changes of  monotonicity and concavity of $u$, the   phase plane analysis of the trajectories associated to the solution $x$ will be repeatedly used. 

A particularly delicate step in the proof of Theorem \ref{main} will be the proof of the non existence of solutions in the critical cases $p=p^*_\pm$, as well as of the uniqueness of the fast decaying solutions, which will be  obtained by using different arguments  for  $\Mp$ and $\Mm$. For $\Mm$, we heavily exploit the fact that $p^*_- >\frac{\tilde{N}_-+2}{\tilde{N}_--2}$ and the proof relies  on some   properties of the solutions of supercritical semilinear problems. For $\Mp$, for which $p^*_+<\frac{\tilde{N}_++2}{\tilde{N}_+-2}$,  a different  proof will be obtained as an application of Gauss-Green Theorem in the phase plane.
\smallskip

The paper is organized as follows: in Section 2 we recall some basic properties of radial solutions of problem \eqref{eq:probMpm} and their Emden-Fowler transform, whereas the existence of positive radial solutions for  $p$ supercritical is proved in Section 3. Section 4 will be then  devoted to the proof of nonexistence of nontrivial solutions when $p$ is strictly subcritical,  while the critical cases $p=p^*_+$ and $p=p^*_-$ are addressed respectively in Section 5 and Section 6. Finally, in Section 7,  we complete the proof of Theorem \ref{main} by showing that problem \eqref{eq:probMpm} has a unique fast decaying solution and infinitely many slow or pseudo-slow decaying solutions, according  to the  initial slope.

\section{Radial solutions of problem \eqref{eq:probMpm} and the Emden-Fowler transform}

In order to study the existence of solutions for problem \eqref{eq:probMpm}, we can assume, without loss of generality, that $B$ is the unit ball of $\R^N$ centered at the origin, by the invariance of the equation with respect to translations and to the scaling  $\tilde u(x):=\gamma^{\frac{2}{p-1}}u(\gamma x)$ for any $\gamma>0$. 
\smallskip

For any $\alpha>0$, let us introduce the initial value problem
\begin{equation}\label{eq:initvalprob+}
\begin{cases}
u^{\prime\prime}(r)=M^+\left(-\frac{\lambda(N-1)}{r}K^+(u^\prime) - u^{p}\right) & \hbox{for} \ r>1,\\
\qquad  \ \ \ \ \ u(r)>0 & \hbox{for} \ r>1,\\
\qquad  \ \ \ \ \ u(1)=0, \ \ u^\prime(1)=\alpha, & 
\end{cases}
\end{equation}
where 
$$M^+(s):=\begin{cases} s/\Lambda & \hbox{if} \ s\geq 0,\\ s/\lambda & \hbox{if} \ s< 0,\end{cases} \ \ \ \ \ K^+(s):=\begin{cases} \frac{\Lambda}{\lambda}s & \hbox{if} \ s\geq 0,\\ s& \hbox{if} \ s< 0.\end{cases}$$
A direct computation shows that $u$ is a radial solution of problem \eqref{eq:probMpm} with ${\mathcal F}=\Mp$ if and only if $u$ satisfies  \eqref{eq:initvalprob+} in $(1,+\infty)$ for some $\alpha>0$.

Analogously, a radial solution of \eqref{eq:probMpm} with ${\mathcal F}=\Mm$ is nothing but a  global solution of 
\begin{equation}\label{eq:initvalprob-}
\begin{cases}
u^{\prime\prime}(r)=M^-\left(-\frac{\Lambda(N-1)}{r}K^-(u^\prime) - u^{p}\right) & \hbox{for} \ r>1,\\
\qquad  \ \ \ \ \ u(r)>0 & \hbox{for} \ r>1,\\
\qquad  \ \ \ \ \ u(1)=0, \ \ u^\prime(1)=\alpha, & 
\end{cases}
\end{equation}
with
$$M^-(s):=\begin{cases} s/\lambda & \hbox{if} \ s\geq 0,\\ s/\Lambda & \hbox{if} \ s< 0,\end{cases} \ \ \ \ \ K^-(s):=\begin{cases} \frac{\lambda}{\Lambda}s & \hbox{if} \ s\geq 0,\\ s& \hbox{if} \ s< 0.\end{cases}$$

For the rest of the section, let us focus on the case ${\mathcal F}=\Mp$. The same observations can be made for 
${\mathcal F}=\Mm$, by just exchanging $\Lambda$ with $\lambda$.

Let us denote by $u_\alpha$ the unique positive maximal solution of \eqref{eq:initvalprob+}, defined on a maximal interval $[1,\rho_\alpha)$, with $\rho_\alpha \leq +\infty$. Some general properties of the function $u_\alpha$ have been established in \cite{GLP}. In particular (see \cite[Lemma 2.1]{GLP}), it is known that, for any $\alpha>0$, there exists a unique $\tau_\alpha\in (1, \rho_\alpha)$ such that $u_\alpha'(r)>0$ for $r\in [1,\tau_\alpha)$ and $u_\alpha'(r)<0$ for $r\in (\tau_\alpha, \rho_\alpha)$. Moreover, if $\rho_\alpha <+\infty$ one has $u_\alpha(\rho_\alpha)=0$, whereas if $\rho_\alpha=+\infty$, then $\lim_{r\to +\infty} u_\alpha(r)=0$ (see \cite[Corollary 2.3]{GLP}). We remark that if $\rho_\alpha<+\infty$, then $u_\alpha$ is a radial solution of the Dirichlet problem
$$
\begin{cases}
-\Mp (D^2 u)=u^{p} & \hbox{in} \ A_{1, \rho_\alpha},\\
\qquad  \ \ \ \ \ u=0 & \hbox{on} \ \partial A_{1, \rho_\alpha},
\end{cases}
$$
where $A_{1, \rho_\alpha}$ is the annular domain with radii $1<\rho_\alpha$, and $u'_\alpha(\rho_\alpha)<0$ by the Hopf boundary maximum principle. On the contrary, if $\rho_\alpha=+\infty$, then $u_\alpha$ is a positive solution in the exterior domain $\R^N\setminus \overline{B_1}$ and, by the results of \cite{AS}, this necessarily implies that $p>\frac{\tilde{N}_+}{\tilde{N}_+-2}$.

We further notice that the function $u_\alpha$,  which satisfies $u''_\alpha(1)<0$,   cannot be globally (or even definitely) concave in $(1,\rho_\alpha)$, since otherwise $\rho_\alpha<+\infty$ and 
$$u^p_\alpha(\rho_\alpha)=-\lambda\left(u''_\alpha(\rho_\alpha)+(N-1)u'_\alpha(\rho_\alpha)/\rho_\alpha\right)>0\, .$$
We will denote by $\sigma_\alpha\in (\tau_\alpha, \rho_\alpha)$ 
the first zero of $u''_\alpha$. We will show below that  $\sigma_\alpha$ actually is the first point where the function $u_\alpha$ changes its concavity.

A useful tool for studying the problem is the so called Emden-Fowler transform (see \cite{F,FQ}), which reduces the initial problem \eqref{eq:initvalprob+} to an autonomous equation. Setting 
\begin{equation}\label{EF}
x(t)=r^{\frac{2}{p-1}}u(r)\, ,\quad r=e^t\, ,
\end{equation}
a direct computation shows that 
$$
\begin{array}{cl}
\displaystyle u'\geq 0 & \displaystyle  \Longleftrightarrow \ x'\geq \frac{2}{p-1} x\, ,\\[2ex]
\displaystyle  u''\leq 0  & \displaystyle  \Longleftrightarrow \ x'\geq \frac{2}{p-1} x-\frac{x^p}{\lambda\, (N-1)}\, .
\end{array}
$$
Hence, the new unknown $x$ satisfies
\begin{equation}\label{Xeq}
\left\{
\begin{array}{ll}
\displaystyle x''= \tilde{a}_- x'+ \tilde{b}_- x-\frac{x^p}{\lambda} & \displaystyle \hbox{ if } \ x'\geq \frac{2}{p-1} x\, , \\[2ex]
\displaystyle x''= a x'+ b x-\frac{x^p}{\lambda} & \displaystyle \hbox{ if } \  \frac{2}{p-1} x\geq x' \geq \frac{2}{p-1} x-\frac{x^p}{\lambda\, (N-1)}\, ,\\[2ex]
\displaystyle x''= \tilde{a}_+ x'+ \tilde{b}_+ x-\frac{x^p}{\Lambda} & \displaystyle \hbox{ if }  \     x'\leq   \frac{2}{p-1} x-\frac{x^p}{\lambda\, (N-1)}\, ,
\end{array}
\right.
\end{equation}
with coefficients defined respectively as
\begin{equation}\label{coeff}
\begin{array}{c}
\displaystyle a : =\frac{N+2 -(N-2)p}{p-1}\, , \qquad \tilde{a}_\pm : =\frac{\tilde{N}_\pm+2 -(\tilde{N}_\pm-2)p}{p-1}\, ,\\[2ex]
\displaystyle b : =\frac{2\left( (N-2)p-N\right)}{(p-1)^2}\, , \qquad \tilde{b}_\pm : =\frac{2\left( (\tilde{N}_\pm-2)p-\tilde{N}_\pm\right)}{(p-1)^2}\, .
\end{array}
\end{equation}

\begin{center}
\includegraphics[scale=0.25]{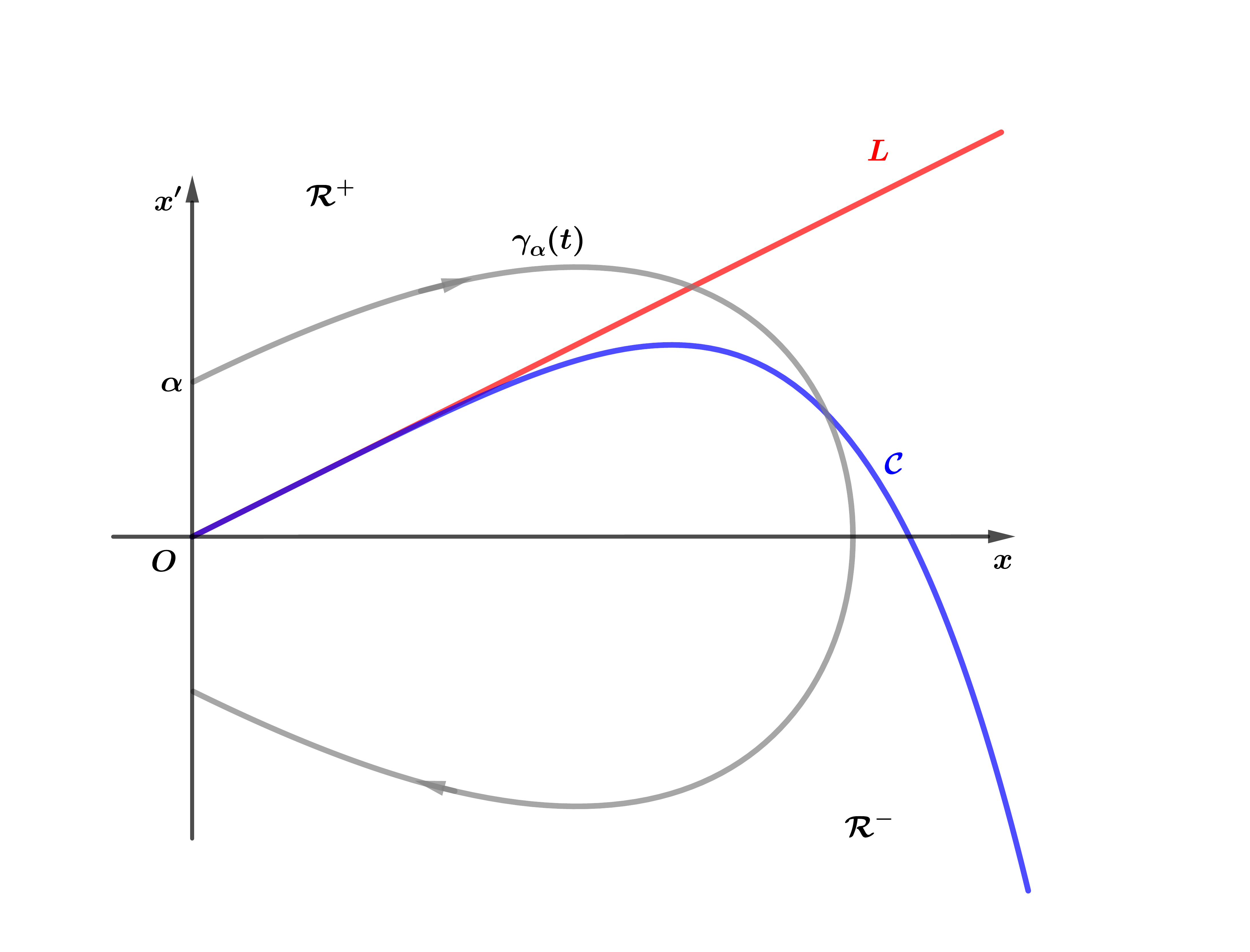}
\end{center}
Associated with a solution $x$ of the above problem, one can consider the trajectory $\gamma (t)=(x(t), x'(t))$ in the phase plane. Thus, the coefficients of the autonomous equation satisfied by $x$ are piecewise constant and  have jumps whenever the trajectory $\gamma$  crosses either the half-line $L=\left\{x'=\frac{2}{p-1}x\right\}$ or the curve $\mathcal{C}=\left\{x'=\frac{2}{p-1}x -\frac{x^p}{\lambda\, (N-1)}\right\}$. We will denote  by $\mathcal{R}^+$ and $\mathcal{R}^-$ the open regions of the right half-plane lying respectively above and below the curve $\mathcal{C}$.

In particular,  the trajectory $\gamma_\alpha(t)=(x_\alpha(t), x'_\alpha(t))$  associated with the transform $x_\alpha$ of the solution $u_\alpha$ is defined on the interval  $[0, \log \rho_\alpha)$, it lays in the right half plane   for $t\in (0, \log \rho_\alpha)$ and  it satisfies the initial condition $\gamma_\alpha(0)=(0,\alpha)$.  We remind the  reader that the function $x_\alpha$ is increasing as long as the trajectory $\gamma_\alpha$ lays in the first quadrant, as well as $x_\alpha$ decreases when $\gamma_\alpha$ is in the fourth quadrant. If $\rho_\alpha<+\infty$, then the trajectory $\gamma_\alpha$ reaches the $x'-$axes at the finite time $t=\log \rho_\alpha$, with $\gamma_\alpha(\log \rho_\alpha)=\left(0, \rho_\alpha^{\frac{p+1}{p-1}}u'_\alpha(\rho_\alpha)\right)$ and  $u_\alpha'(\rho_\alpha)<0$. The typical behavior of a trajectory $\gamma_\alpha$ having $\rho_\alpha<+\infty$ is depicted in the figure above.

In the following result we collect some properties of the trajectory $\gamma_\alpha$, partially observed in previous contributions (we refer in particular to \cite{FQ, BGLP, GLP}).  

 \begin{lemma}\label{gammaalfa}
 For each $\alpha>0$, there exist unique $0<t_\alpha<s_\alpha<\log \rho_\alpha$ such that
 
 \begin{itemize}
 \item[(i)] $x'_\alpha (t)>\frac{2}{p-1}x_\alpha(t)$ for all $t\in [0, t_\alpha)$ and $x'_\alpha (t)<\frac{2}{p-1}x_\alpha(t)$ for all $t\in (t_\alpha, \log \rho_\alpha)$;
 \smallskip
 
 \item[(ii)] $\frac{2}{p-1}x_\alpha(t) -\frac{x_\alpha^p(t)}{\lambda\, (N-1)}< x'_\alpha (t)<\frac{2}{p-1}x_\alpha(t)$ for all $t\in (t_\alpha, s_\alpha)$ and $\gamma_\alpha(s_\alpha)\in \mathcal{C}$;
 \smallskip
 
 \item[(iii)]  the trajectory $\gamma_\alpha$, as any other trajectory $\gamma$, can intersect the curve $\mathcal{C}$ from above only transversally and only at points satisfying $x^{p-1}> \frac{\lambda\, (N-1)}{p}$. In particular,  $x^{p-1}_\alpha(s_\alpha)> \frac{\lambda\, (N-1)}{p}$ and $\gamma_\alpha(t)\in \mathcal{R}^-$  for $t$ in a right neighborhood of $s_\alpha$;
 \smallskip
 
 \item[(iv)]the trajectory $\gamma_\alpha$,  as any other trajectory $\gamma$, can intersect the curve $\mathcal{C}$ from below only at points satisfying $x^{p-1}\leq \frac{\lambda\, (N-1)}{p}$ and intersections from below are always transversal if they occur at points satisfying  $x^{p-1} < \frac{\lambda\, (N-1)}{p}$;
 \smallskip
 
 \item[(v)] if either $\rho_\alpha<+\infty$ or $\rho_\alpha=+\infty$ and $\gamma_\alpha(t)\to (0,0)$ as $t\to +\infty$, then  the trajectory $\gamma_\alpha$ intersects (and crosses) the $x$-axis exactly once and $\gamma_\alpha(t)\in \mathcal{R}^-$ for all $t\in (s_\alpha , \log \rho_\alpha )$;
 \smallskip
 
 \item[(vi)]  the map  $\gamma_\alpha(t)$ is bounded for $t\in [0, \log \rho_\alpha)$;
 \smallskip
 
 \item [(vii)] if $\rho_\alpha=+\infty$, then, as $t\to +\infty$, only three possible cases can occur: either $\gamma_\alpha(t)\to (0,0)$, or $\gamma_\alpha(t)\to \left( (\tilde{b}_+\Lambda)^{1/(p-1)},0\right)$, or $\gamma_\alpha(t)$ approaches the  trajectory of a periodic orbit lying in the right half-plane and enclosing the equilibrium point $\left( (\tilde{b}_+\Lambda)^{1/(p-1)},0\right)$;
\smallskip

\item[(viii)] if $\rho_\alpha=+\infty$ and $\gamma_\alpha(t)\to (0,0)$ as $t\to +\infty$, then there exist  non zero constants $c_{1,\alpha}$ and $c_{2,\alpha}$ such that
$$
e^{-\lambda_1 \, t} \gamma_\alpha (t)\to (c_{1,\alpha},c_{2,\alpha}) \quad \hbox{ as } t\to +\infty
$$
where $\lambda_1 =-\frac{ (\tilde{N}_+-2)p-\tilde{N}_+}{p-1}$;
\smallskip

\item[(ix)] if $\rho_\alpha=+\infty$ and $\gamma_\alpha(t)$ approaches the  trajectory of a periodic orbit as $t\to +\infty$, then $\gamma_\alpha(t)$ intersects the curve $\mathcal{C}$ infinitely many times for $t\in [0,+\infty)$.

 \end{itemize}
 \end{lemma}

\proof By setting $t_\alpha=\log\tau_\alpha$ and $s_\alpha=\log \sigma_\alpha$, properties (i) e (ii) immediately follow by the definition of $\tau_\alpha$ and $\sigma_\alpha$.

In order to prove  (iii), let us assume that, for a time  $s\in (0, \log \rho_\alpha)$,  one has $\gamma_\alpha(s)\in \mathcal{C}$ and $\gamma_\alpha (t)\in \mathcal{R}^+$ for $t$ in a left neighborhood of $s$. This means that
 $$
 x_\alpha'(s)=\frac{2}{p-1}x_\alpha(s) -\frac{x_\alpha^p(s)}{\lambda\, (N-1)}
 $$
 and 
$$
x_\alpha'(t)>\frac{2}{p-1}x_\alpha(t) -\frac{x_\alpha^p(t)}{\lambda\, (N-1)} \quad \hbox{ for } t\in (s-\delta, s)\, ,
$$
for some $\delta>0$. It then follows that
$$
x_\alpha''(s)\leq \left( \frac{2}{p-1} -\frac{p \, x_\alpha^{p-1}(s)}{\lambda\, (N-1)}  \right) x'_\alpha(s)
= \left( \frac{2}{p-1} -\frac{p \, x_\alpha^{p-1}(s)}{\lambda\, (N-1)}\right)  \left( \frac{2}{p-1} -\frac{x_\alpha^{p-1}(s)}{\lambda\, (N-1)}
\right) x_\alpha(s)\, .
$$
On the other hand, by \eqref{Xeq} one has
$$
x_\alpha''(s)= a \, x'_\alpha(s)+b\,  x_\alpha(s)-\frac{x^p_\alpha(s)}{\lambda} = \left( \frac{4}{(p-1)^2}-\frac{p+3}{\lambda \, (N-1)\, (p-1)} x_\alpha^{p-1}(s)\right) x_\alpha(s)\, .
$$
 Combining the above relationships, since $x_\alpha(s)>0$, we then obtain
 $$
 x_\alpha^{p-1}(s)\geq \frac{\lambda\, (N-1)}{p}\, ,
 $$
 and the intersection is tangential if and only if $ x_\alpha^{p-1}(s)= \frac{\lambda\, (N-1)}{p}$. 
But a second order analysis at the intersection point, as in the proof of Lemma 3.1 in \cite{FQ}, shows that intersections from above  never occur at the point $\left( \left(\frac{\lambda\, (N-1)}{p}\right)^{\frac{1}{p-1}}, \frac{p+1}{p\, (p-1)} \left(\frac{\lambda\, (N-1)}{p}\right)^{\frac{1}{p-1}}\right) $. This proves statement (iii). The same argument, with reversed inequalities, proves claim (iv).

In order to prove (v), let $\gamma_\alpha$ be such that either $\rho_\alpha<+\infty$ or $\gamma_\alpha(t)\to (0,0)$ as $t\to +\infty$. Since $\gamma_\alpha(0)=(0,\alpha)$, then the trajectory crosses at least once  the $x$-axis. Let us call $t_0\in (0,\log \rho_\alpha )$ the  time of the first intersection between $\gamma_\alpha$ and the $x$-axis, that is  $x'_\alpha(t_0)=0$ and $x'_\alpha(t)>0$ for all $t\in [0, t_0)$. Then, $x_\alpha''(t_0)\leq0$. On the other hand, by \eqref{Xeq}, one has
$$
x_\alpha''(t_0)=\left\{
\begin{array}{ll}
x_\alpha(t_0)\left( b-\frac{x_\alpha(t_0)^{p-1}}{\lambda}\right) & \hbox{if } x_\alpha(t_0)^{p-1}\geq \frac{2 \lambda (N-1)}{p-1}\, ,\\[1ex]
x_\alpha(t_0)\left( \tilde b_+ -\frac{x_\alpha(t_0)^{p-1}}{\Lambda}\right) & \hbox{if } x_\alpha(t_0)^{p-1}< \frac{2 \lambda (N-1)}{p-1}\, .
\end{array}
\right.
$$
This immediately implies that $x_\alpha''(t_0)<0$ if $x_\alpha(t_0)^{p-1}\geq \frac{2 \lambda (N-1)}{p-1}$, since $\frac{2 \lambda (N-1)}{p-1}>\lambda \, b$, whereas it yields $x_\alpha (t_0)^{p-1} \geq \Lambda \tilde b_+$ in the case  $x_\alpha(t_0)^{p-1}< \frac{2 \lambda (N-1)}{p-1}$. But if  $x_\alpha(t_0)^{p-1} = \Lambda \tilde{b}_+$,  then $\tilde{b}_+>0$ and $\gamma_\alpha(t_0)=\left( (\Lambda \tilde{b}_+)^{\frac{1}{p-1}}, 0\right)$, which is impossible since   $\left( (\Lambda \tilde{b}_+)^{\frac{1}{p-1}}, 0\right)$ is an equilibrium point. Hence, in both cases we obtain $x_\alpha''(t_0)<0$ and, for $t>t_0$, the trajectory enters the fourth quadrant and $x_\alpha(t)$ starts to decrease. We then observe that for $t>t_0$ the trajectory   never intersects  the $x$-axis, since otherwise, arguing as above,  it would cross it, entering again into the first quadrant, and, in order to avoid self-intersection, it could not reach any point of the $x'$-axes for subsequent times. This proves that $\gamma_\alpha$ intersects exactly once the $x$- axis. We notice further that, if $s_\alpha <t_0$, that is if 
 the  first intersection point between $\gamma_\alpha$ and $\mathcal{C}$ lies in the first quadrant,  then, for $t\in (s_\alpha,t_0]$ the trajectory lies in the first quadrant and  $x_\alpha(t)>x_\alpha(s_\alpha)$, so that   $\gamma_\alpha(t)$ cannot intersect $\mathcal{C}$ by property (iv). Therefore,  $\gamma_\alpha(t)$ stays in $\mathcal{R}^-$ until it reaches the $x$-axis at the time $t_0$, then  it  enters the fourth quadrant,  $x_\alpha(t)$ becomes decreasing, and $\gamma_\alpha(t)\in \mathcal{R}^-$ also for $t\geq t_0$. The same argument can be applied in the case $s_\alpha\geq t_0$,  since then $x'_\alpha(t)<0$ for any $t\in (s_\alpha, \log \rho_\alpha)$ and, by  property (iv), $\gamma_\alpha$ cannot intersect $\mathcal{C}$ from below.  This completely proves the claim.

Let us now prove (vi). The statement is trivial if $\rho_\alpha<+\infty$, so we assume $\rho_\alpha=+\infty$ and, thus, $p>\frac{\tilde{N}_+}{\tilde{N}_+-2}$.

 From (iv) it follows that, if  $\gamma_\alpha(t)$ intersects again the curve $\mathcal{C}$ for $t>s_\alpha$, and if we call
 $\hat{s} >s_\alpha $ the first time such that $\gamma_\alpha(\hat s)\in \mathcal{C}$ and $\gamma_\alpha(t)\in \mathcal{R}^-$ for $t\in (s_\alpha, \hat s)$, then necessarily $x'_\alpha(\hat s)>0$. Hence, $\gamma_\alpha(t)$ intersects twice the $x$-axis for $t\in [0,\hat s]$ and therefore, for all $t\in [0,+\infty)$, one has
 $$|\gamma_\alpha(t)|\leq \max_{s\in [0,\hat s]}|\gamma_\alpha(s)|\, .
 $$
 It remains to prove the boundedness of the trajectories $\gamma_\alpha (t)$ satisfying $\gamma_\alpha(t)\in \mathcal{R}^-$ for all $t>s_\alpha$. In this case it is convenient to look back at the function $u_\alpha$,  which is convex and decreasing in $(\sigma_\alpha,+\infty)$. By  \eqref{eq:initvalprob+}, $u_\alpha$ satisfies
 $$
 u''_\alpha +\frac{\tilde{N}_+-1}{r} u'_\alpha =-\frac{u_\alpha^p}{\Lambda}\quad \hbox{ for } r\in [\sigma_\alpha,+\infty)\, ,
 $$
 which can be written as
 \begin{equation}\label{div1}
 \left( r^{\tilde{N}_+-1} u'_\alpha\right)'= -r^{\tilde{N}_+-1} \frac{u_\alpha^p}{\Lambda}\quad \hbox{ in }  [\sigma_\alpha,+\infty)\, .
 \end{equation}
Integrating between $\sigma_\alpha$ and any $r>\sigma_\alpha$, and using the decreasing monotonicity of $u_\alpha$, we get
$$
r^{\tilde{N}_+-1}u'_\alpha(r)\leq \int_{\sigma_\alpha}^r\left( s^{\tilde{N}_+-1} u'_\alpha(s)\right)' ds\leq - \frac{u_\alpha^p(r)}{\Lambda\, \tilde{N}_+}\left( r^{\tilde{N}_+}-\sigma_\alpha^{\tilde{N}_+}\right)\, .
$$
By integrating once again, we deduce
$$
-\frac{u_\alpha^{1-p}(r)}{p-1}\leq \int_{\sigma_\alpha}^r\frac{u'_\alpha(s)}{u^p_\alpha(s)} ds\leq -\frac{1}{\Lambda\, \tilde{N}_+}
 \int_{\sigma_\alpha}^r \frac{s^{\tilde{N}_+}-\sigma_\alpha^{\tilde{N}_+}}{s^{\tilde{N}_+-1}} ds\leq - \frac{1}{2 \Lambda\, \tilde{N}_+}
 \left( r^2-\frac{\tilde{N}^+}{\tilde{N}^+-2}\sigma_\alpha^2\right)\, .
 $$
 Hence, for $r\geq \hat{\sigma}_\alpha: = \sqrt{\frac{2\tilde{N}^+}{\tilde{N}^+-2}} \sigma_\alpha$, we have
 \begin{equation}\label{est1}
 u^{p-1}_\alpha (r)\leq \frac{4 \Lambda\, \tilde{N}_+}{(p-1)\, r^2}=: \frac{c_1}{r^2}\, .
 \end{equation}
 Using \eqref{est1} into \eqref{div1}, we also have
 $$
 \left( r^{\tilde{N}_+-1} u'_\alpha\right)' \geq -\frac{c_1^{\frac{p}{p-1}}}{\Lambda} r^{\tilde{N}_+ -\frac{2p}{p-1}-1}\quad  \hbox{ in }  [\hat{\sigma}_\alpha,+\infty)\, .
 $$
By integrating in $[\hat{\sigma}_\alpha, r]$ and recalling that $p>\frac{\tilde{N}_+}{\tilde{N}_+-2}$, we get
$$
r^{\tilde{N}_+-1} u'_\alpha(r)\geq \hat{\sigma}_\alpha^{\tilde{N}_+-1} u'_\alpha(\hat{\sigma}_\alpha)- 
\frac{(p-1)c_1^{\frac{p}{p-1}}}{\Lambda\, \left( p(\tilde{N}_+-2)-\tilde{N}_+\right)} r^{\frac{p(\tilde{N}_+-2)-\tilde{N}_+}{p-1}} \, .
$$
Therefore, for $r\geq  \hat{\sigma}_\alpha$, using again that $p>\frac{\tilde{N}_+}{\tilde{N}_+-2}$, we obtain the estimate
\begin{equation}\label{est2}
0\leq -u'_\alpha(r) \leq \left( -  \hat{\sigma}_\alpha^{\tilde{N}_+-1} u'_\alpha(\hat{\sigma}_\alpha)+ \frac{(p-1)c_1^{\frac{p}{p-1}}}{\Lambda\, \left( p(\tilde{N}_+-2)-\tilde{N}_+\right)} \right) \frac{1}{r^\frac{p+1}{p-1}}=: \frac{c_2}{r^\frac{p+1}{p-1}}\, .
\end{equation}
Recalling the definition \eqref{EF} of $x_\alpha$, estimates \eqref{est1} and \eqref{est2} then yield
$$
|\gamma_\alpha(t)| \leq c_3 \quad \hbox{ for } t\geq \log \hat{\sigma}_\alpha\ ,
$$
for a positive constant $c_3$ depending on $p, \Lambda, \lambda, N$ and $\alpha$. This proves (vi).

Claim (vii) then follows by (vi) and the classical Poincar\'e-Bendixson Theorem (see e.g. \cite{H}), after observing that the only equilibrium points of system \eqref{Xeq} reachable by trajectories lying in the right half-plane are $(0,0)$ and $\left( (\tilde{b}_+\Lambda)^{1/(p-1)},0\right)$ if $p>\frac{\tilde{N}_+}{\tilde{N}_+-2}$ (i.e. $\tilde{b}_+>0$), which is always the case if $\rho_\alpha=+\infty$.

Also statement (viii) classically follows from the perturbed linear systems theory, $\lambda_1$ being the negative eigenvalue of the linear system obtained by linearizing around zero   system \eqref{Xeq}. For a detailed proof, we refer e.g. to Lemma 3.3 in \cite{K}.

Let us finally prove (ix). 
We first claim that if   a periodic orbit $\gamma$ of system \eqref{Xeq}  is the $\omega$-limit set of a trajectory $\gamma_\alpha$, then $\gamma$ intersects the curve $\mathcal{C}$. Indeed, the trajectory of such a periodic orbit, if any,   is a closed curve $\gamma$ contained in the right half-plane and winding around the equilibrium point $\left( (\tilde{b}_+\Lambda)^{1/(p-1)},0\right)$. If, by contradiction, $\gamma$ does not intersect $\mathcal{C}$, then $\gamma (t)$ lies entirely in the region $\mathcal{R}^-$. This implies that $\gamma$ is the trajectory associated with a positive periodic solution $x(t)$ of the single equation
$$
x''= \tilde{a}_+ x' +\tilde{b}_+ x -\frac{x^p}{\Lambda}\, .
$$
 Hence, the energy function
\begin{equation}\label{energy}
E(x,y)=\frac{y^2}{2}-\tilde{b}_+ \frac{x^2}{2}+\frac{x^{p+1}}{\Lambda\, (p+1)}\, ,
\end{equation}
when evaluated along   $\gamma(t)$, satisfies
$$
\frac{d}{dt} \left( E(\gamma(t))\right) = \tilde{a}_+ (x')^2\, .
$$
In particular, if $\tilde{a}_+\neq 0$, i.e. if $p\neq  \frac{\tilde{N}_++2}{\tilde{N}_+-2}$, then the  function $E(\gamma(t))$ is strictly monotone and periodic: a contradiction. On the other hand, if $\tilde{a}_+= 0$, i.e. if $p= \frac{\tilde{N}_++2}{\tilde{N}_+-2}$, then $E(\gamma(t))$ is constant, that is $\gamma$ is a level line of the function $E$ contained in $\mathcal{R}^-$. Every level line of the function $E$ contained in $\mathcal{R}^-$ actually is a periodic orbit of system \eqref{Xeq}, and,  since $\gamma$ is contained in $\mathcal{R}^-$, by continuity there exists a larger level line of the function $E$ still contained in $\mathcal{R}^-$ and enclosing $\gamma$. This means that $\gamma$ cannot be approached by any of the trajectories $\gamma_\alpha$, again a contradiction to the assumptions. 

Therefore, for any value of $p>\frac{\tilde{N}_+}{\tilde{N}_+-2}$, the $\omega$-limit $\gamma$ intersects  $\mathcal{C}$ at least once.  Hence, if  $\gamma_\alpha $ is a trajectory winding around $\gamma $ as $t\to +\infty$,  then $\gamma_\alpha (t)$ intersects $\mathcal{C}$ infinitely many times as $t$ ranges in $[0,+\infty)$ (see also the proof of Proposition 3.3 in \cite{FQ}).

\hfill$\Box$

By means of Lemma \ref{gammaalfa} (vii), the maximal solutions $u_\alpha$ of problem \eqref{eq:initvalprob+} for which $\rho_\alpha=+\infty$ can be classified in terms of the asymptotic behavior of the associated trajectory $\gamma_\alpha$. Precisely, a solution $u_\alpha$ is called fast, slow or pseudo-slow decaying as $r\to +\infty$ provided that $\gamma_\alpha(t)$ converges as $t\to +\infty$ respectively to $(0,0)$, $\left( (\tilde{b}_+\Lambda)^{1/(p-1)},0\right)$ or  to a  periodic orbit.
Furthermore, thanks to Lemma \ref{gammaalfa} (viii), a fast decaying solution may be equivalently defined as a solution $u_\alpha : [1,+\infty)\to \R$ satisfying the asymptotic condition \eqref{fast}.

\begin{remark}\label{entire1}
{\rm Positive entire radial solutions of equations 
$$
-\Mpm (D^2u)= u^p\quad \hbox{ in } \R^N
$$
have been studied in detail in \cite{FQ}. In this case the initial value problems for the ODE to be considered are
$$
\begin{cases}
v^{\prime\prime}(r)=M^\pm\left(-\frac{\lambda(N-1)}{r}K^\pm(v^\prime) - v^{p}\right) & \hbox{for} \ r>0,\\
\qquad  \ \ \ \ \ v(r)>0 & \hbox{for} \ r>0,\\
\qquad  \ \ \ \ \ v(0)=\alpha, \ \ v^\prime(0)=0, & 
\end{cases}
$$
and for the maximal solutions $v_\alpha$,   defined on $[0,r_\alpha)$ for some $0<r_\alpha \leq +\infty$, one has $r_\alpha<+\infty$ if and only if $p<p^\ast_\pm$. 
We recall that, by the homogeneity property of the equation, for any $\alpha, \alpha'>0$ one has 
\begin{equation}\label{scaling}
v_{\alpha'}(r)=\frac{\alpha'}{\alpha} v_\alpha \left(\left( \frac{\alpha'}{\alpha}\right)^{\frac{p-1}{2}}r\right)\, ,
\end{equation}
that is changing the initial condition just reflects in scaling the solution. Therefore, the property $r_\alpha<+\infty$ or $r_\alpha=+\infty$ only depends on $p$.

The trajectory $\Gamma_\alpha (t)=\left(\xi_\alpha (t), \xi_\alpha '(t)\right)$ associated with the Emden-Fowler transform $\xi_\alpha (t)$ of a solution $v_\alpha(r)$, defined through  \eqref{EF}, is defined for $t\in (-\infty , \log r_\alpha)$ and satisfies
$$
\lim_{t\to -\infty} \Gamma_\alpha(t)= (0,0)\, .
$$
Relationship \eqref{scaling} implies that
$$
\Gamma_{\alpha'}(t)=\Gamma_\alpha \left( t+\frac{p-1}{2} \log \frac{\alpha'}{\alpha}\right)\, ,
$$
and thus the initial value $\alpha$ for the solution $v_\alpha(r)$ does not affect the support of the trajectory $\Gamma_\alpha(t)$, which is then denoted simply by $\Gamma (t)$. 
In the phase plane, the trajectory $\Gamma  (t)$ is  a curve exiting from the origin, lying  always below the line $L$, and staying above the curve $\mathcal{C}$ for $t\in (-\infty, S_\alpha)$, for some $S_\alpha< \log r_\alpha$. 
We  recall from \cite{FQ} that
$\Gamma (t)$ satisfies properties (iii)-(ix) of Lemma \ref{gammaalfa}, and fast, slow or pseudo-slow decaying solutions $v_\alpha$ are accordingly defined. A crucial result proved in \cite{FQ} is that, independently of $\alpha$, 
 $v_\alpha$ is a fast decaying solution if and only if $p=p^\ast_\pm$. 
}
\end{remark}

\section{Existence in the supercritical case $p>p^*_\pm$ for $\Mpm$}

This section is devoted to the proof of the existence statement in Theorem \ref{main}.

\begin{theorem}\label{mainteo:Existence}
If $\F=\Mp$ and $p>p^*_+$ then \eqref{eq:probMpm} has a nontrivial radial solution.\\
If $\F=\Mm$ and $p>p^*_-$ then \eqref{eq:probMpm} has a nontrivial radial solution.
\end{theorem}
\begin{proof} 
We write the proof for  $\Mp$. Unless otherwise said, the proof for $\Mm$ is obtained just by exchanging  $\lambda$ with $\Lambda$. 

Assume by contradiction that the thesis is false. Then, for any $\alpha>0$,  the unique maximal solution $u_\alpha=u_\alpha(r)$ of the initial value problem \eqref{eq:initvalprob+} is defined on the finite interval $[1, \rho_\alpha]$. We recall that 
 there exists  a unique $\tau_\alpha \in (1,\rho_\alpha)$ such that $u_\alpha^\prime(\tau_\alpha)=0$, with $u_\alpha^\prime(r)>0$ for $r\in [1,\tau_\alpha)$ and   $u_\alpha^\prime(r)<0$ for $r\in (\tau_\alpha,\rho_\alpha]$. Moreover, by Lemma \ref{gammaalfa} (ii)-(iii)-(v), there exists a unique $\sigma_\alpha\in (\tau_\alpha ,\rho_\alpha)$ such that $u_\alpha''(\sigma_\alpha)=0$, with $u_\alpha''(r)<0$ for $r\in [1,\sigma_\alpha)$ and   $u_\alpha''>0$ for $r\in (\sigma_\alpha,\rho_\alpha]$.
\smallskip

Next, we perform the asymptotic analysis as $\alpha \to 0^+$.
 From \cite[Lemma 3.1]{GLP} we know that $\tau_\alpha\to +\infty$, $u_\alpha(\tau_\alpha)\to 0$.  Hence, we have also $\sigma_\alpha \to +\infty$, $\rho_\alpha \to +\infty$ as $\alpha\to 0^+$. Setting $m_\alpha:=u_\alpha(\tau_\alpha)$, we consider the   rescaled function
$$\tilde u_\alpha(r):=\frac{1}{m_\alpha} u_\alpha\left(\frac{r}{m_\alpha^{\frac{p-1}{2}}}\right), \ \ r\in\left[m_\alpha^{\frac{p-1}{2}}, m_\alpha^{\frac{p-1}{2}}\rho_\alpha \right].$$
For notational convenience, let us  set  $\tilde \tau_\alpha:= m_\alpha^{\frac{p-1}{2}} \tau_\alpha$, $\tilde \sigma_\alpha:= m_\alpha^{\frac{p-1}{2}}  \sigma_\alpha$, $\tilde \rho_\alpha:=m_\alpha^{\frac{p-1}{2}}  \rho_\alpha$.
Note that, by construction,  $\tilde u_\alpha$ is a radial solution to
\begin{equation}\label{eq:probrescal}
\begin{cases}
-\Mp (D^2 u)=u^{p} & \hbox{in} \ \tilde A_\alpha,\\
\qquad  \ \ \ \ \ u>0 & \hbox{in} \ \tilde A_\alpha,\\
\qquad  \ \ \ \ \ u=0 & \hbox{on} \ \partial \tilde A_\alpha,
\end{cases}
\end{equation}
where $\tilde A_\alpha=\left\{x\in\mathbb R^N:\ m_\alpha^{\frac{p-1}{2}}<|x|< \tilde{\rho}_\alpha\right\}$, and one has $1=\tilde{u}_\alpha(\tilde{\tau}_\alpha)=\max_{\tilde{A}_\alpha} \tilde{u}_\alpha$.
\smallskip

\noindent\textbf{Step 1}: %Up to a subsequence, 
As $\alpha\to 0^+$, we have $\tilde \tau_\alpha \to 0^+$.

Let us consider the functional
$ E_1:[1,\tau_\alpha]\to\R$ defined by
$$ E_1(r):=\frac{(u_{\alpha}^\prime(r))^2}{2}+\frac{u_{\alpha}^{p+1}(r)}{\lambda(p+1)}.$$
A direct computation, as in the proof of \cite[Lemma 3.1]{GLP}, shows that
 $E_1^\prime(r)\leq 0$ for $ r \in [1,\tau_\alpha]$.\\
Thus, from $E_1(r)\geq E_1(\tau_\alpha)$ and $u_\alpha'(r)\geq 0$ for $r\in [1, \tau_\alpha]$, we infer that
$$u_\alpha^\prime(r)\geq \sqrt{\frac{2}{\lambda(p+1)}\left(m_\alpha^{p+1}-u_\alpha^{p+1}(r)\right)}\qquad\forall \, r\in[1,\tau_\alpha].$$
Integrating between $1$ and $\tau_\alpha$ we get that
\beq\label{eq1:teo1}
\displaystyle\int_1^{\tau_\alpha} \frac{u_\alpha^\prime(r)}{ \sqrt{m_\alpha^{p+1}-u_\alpha^{p+1}(r)}} \ dr \geq  \sqrt{\frac{2}{\lambda(p+1)}}\left(\tau_\alpha-1\right).
\eeq
Performing the change of variable $t=\frac{u_\alpha(r)}{m_\alpha}$, we rewrite \eqref{eq1:teo1} as
\beq\label{eq2:teo1}
\frac{1}{m_\alpha^{\frac{p-1}{2}}}\int_0^{1} \frac{1}{ \sqrt{1-t^{p+1}}} \ dt \geq  \sqrt{\frac{2}{\lambda(p+1)}}\left(\tau_\alpha-1\right).
\eeq
Since $\int_0^{1} \frac{1}{ \sqrt{1-t^{p+1}}}\ dt \leq \int_0^{1} \frac{1}{ \sqrt{1-t^{2}}}\ dt= \frac{\pi}{2}$ and $m_\alpha\to 0$ as $\alpha\to 0^+$, from  \eqref{eq2:teo1} we deduce that  
$$\limsup_{\alpha\to0^+}\tilde\tau_\alpha\leq \frac{\pi}{2}  \sqrt{\frac{\lambda(p+1)}{2}}.$$
Assume now, by contradiction, that $l:=\limsup_{\alpha\to0^+}\tilde\tau_\alpha>0$ and let us select a sequence, still denoted by $\alpha$, such that  $\tilde\tau_\alpha \to l$ as $\alpha\to 0^+$. 

By \eqref{eq:initvalprob+} and the concavity/monotonicity properties of $\tilde{u}_\alpha$, it follows that $\tilde{u}_\alpha$ in particular satisfies
$$
 \tilde u_\alpha^{\prime\prime} +  \frac{\tilde{N}_- -1}{r} \tilde u_\alpha^{\prime}=-\frac{\tilde u_\alpha^{p}}{\lambda} \quad \hbox{ in } 
[m_\alpha^{\frac{p-1}{2}}, \tilde\tau_\alpha]\, .
$$
Since $0\leq \tilde{u}_\alpha \leq 1$, it follows that, up to a further subsequence if necessary, $\tilde u_\alpha(r) \to \tilde u(r)$ in $C^2_{loc}((0,l))$, for some $\tilde u$ satisfying $0\leq \tilde u \leq 1$, $\tilde{u}'\geq 0$ in $(0, l]$,  $\tilde{u}(l)=1$ and
$$
 \tilde u^{\prime\prime} +  \frac{\tilde{N}_- -1}{r} \tilde u^{\prime}=-\frac{\tilde u^{p}}{\lambda} \quad \hbox{ in } 
(0,l]\, .
$$
We observe that $\tilde{u}$ cannot be identically equal to $1$,    since the constant function $1$ does not satisfy the above equation. 
Hence, there exists $r_0\in (0,l)$ such that $\tilde{u}'(r_0)>0$. Observing further that
$$ (\tilde u^\prime r^{\Nm-1})^\prime=-\frac{1}{\lambda} \tilde u^pr^{\Nm-1} \leq 0\, ,$$
we infer that, for all $r\in (0,r_0]$,
$$\tilde u^\prime(r) \geq  \tilde u^\prime(r_0) \frac{r_0^{\Nm-1}}{r^{\Nm-1}}\, .$$
For any $\varepsilon \in (0,r_0)$, let us  integrate  the above inequality from $\varepsilon$ to $r\in \left [ \varepsilon, r_0 \right]$. We obtain
$$
1\geq \tilde u(r) - \tilde u(\e) \geq  \frac{\tilde u^\prime(r_0) r_0^{\Nm-1}}{\Nm-2} \left(\frac{1}{\e^{\Nm-2}}-\frac{1}{r^{\Nm-2}}\right),
$$
which yields  a contradiction  as $\e \to 0^+$. The proof of Step 1 is complete.\\
\smallskip

\noindent\textbf{Step 2}: Up to a  sequence $\alpha \to 0^+$, we have $\tilde \sigma_\alpha \to l$  for some $l \in (0,+\infty)$.

We first observe that, by the definition of $\sigma_\alpha$ and of $x_\alpha$ given by \eqref{EF}, from Lemma \ref{gammaalfa}-(iii) it follows that
$$
\sigma_\alpha^2 u_\alpha^{p-1}(\sigma_\alpha)> \frac{\lambda\, (N-1)}{p}\, .
$$
Recalling the definition of $\tilde{\sigma}_\alpha$ and observing that $u_\alpha(\sigma_\alpha)\leq m_\alpha$, this readily implies that, for every $\alpha>0$,
\begin{equation}\label{eqclaim1teo1}
\tilde{\sigma}_\alpha > \sqrt{\frac{\lambda\, (N-1)}{p}}\, .
\end{equation}
Therefore, in order to prove Step 2, it is enough to show that 
  $\limsup_{\alpha\to 0^+} \tilde \sigma_\alpha < +\infty$.

Assume by contradiction that for some sequence $\alpha \to 0^+$ one has  $\tilde \sigma_\alpha \to +\infty$. Since $\tilde u_\alpha$ is positive and concave in $[\tilde \tau_\alpha, \tilde \sigma_\alpha]$, from $\tilde u_\alpha(\tilde\tau_\alpha)=1$ we infer that
\beq\label{eq5:teo1}
\tilde u_\alpha (r) \geq \frac{\tilde\sigma_\alpha - r}{\tilde\sigma_\alpha -\tilde \tau_\alpha} \ \ \hbox{for any} \ r\in [\tilde \tau_\alpha, \tilde \sigma_\alpha].
\eeq
On the other hand,  in the present assumption  the limit domain of $\tilde u_\alpha$ is $(0,+\infty)$, and, by standard elliptic estimates, up to a further subsequence, $\tilde u_\alpha \to \tilde u$ in $C^2_{loc}((0,+\infty))$, where $\tilde u$ is a radial solution of
\begin{equation}\label{eq:probLimMpm}
\begin{cases}
-\Mp (D^2 w)=w^{p} & \hbox{in} \ \R^N\setminus \{0\},\\
\qquad  \ \ \ \ \ w\geq 0 & \hbox{in} \ \R^N\setminus \{0\}.
\end{cases}
\end{equation}
Passing to the limit as $\alpha\to 0^+$ in \eqref{eq5:teo1} and using Step 1, we deduce that
$$ \tilde u(r)\geq 1 \ \ \hbox{for any} \ r \in (0,+\infty).$$
But $\tilde u \leq 1$ by construction,  and, therefore, $\tilde u \equiv 1$, which is absurd because $\tilde u$ solves \eqref{eq:probLimMpm}. The proof of Step 2 is complete.\\
\smallskip

\noindent\textbf{Step 3}: %Up to a subsequence, 
As $\alpha\to 0^+$, we have $\tilde \rho_\alpha \to +\infty$.

Let us assume, by contradiction, that $k:=\liminf_{\alpha\to 0^+} \tilde{\rho}_\alpha <+\infty$. Then, there exists a sequence still denoted by $\alpha\to 0^+$ such that 
$\tilde \rho_\alpha \to k$. By  \eqref{eqclaim1teo1}, we have that $k>0$.  Recall that, by Step 1, $\tilde{\tau}_\alpha\to 0$ and, possibly considering a further subsequence, $\tilde{\sigma}_\alpha\to l_1$ for some $0<l_1\leq k$. 
\\ By construction and standard elliptic estimates,  it follows that, possibly up to a further subsequence,  $\tilde u_\alpha \to \tilde u$ in $C^2_{loc} ((0,k))$, where $\tilde u$ is in this case  a radial solution of
\begin{equation}\label{eq:problimitMpm}
\begin{cases}
-\Mp (D^2 w)=w^{p} & \hbox{in} \ B_k\setminus\{0\},\\
\qquad  \ \ \ \ \ w\geq 0 & \hbox{in} \  B_k\setminus\{0\},\\
\qquad  \ \ \ \ \ w=0 & \hbox{on} \ \partial B_k.
\end{cases}
\end{equation}
We claim that $\tilde u$ can be extended to a smooth non-trivial radial solution of \eqref{eq:problimitMpm} in the whole ball $B_k$. 
Indeed, by  \eqref{eq:initvalprob+},  in $[\tilde\tau_\alpha,\tilde \sigma_\alpha]$  the function $\tilde u_\alpha$ satisfies 
$$
\tilde u_\alpha^{\prime\prime}(r)+\frac{N-1}{r}\tilde u_\alpha^{\prime}(r)=-\frac{\tilde u_\alpha^p(r)}{\lambda}\,,
$$
and, therefore,
\begin{equation}\label{eq8}
\left(r^{N-1}\tilde u_\alpha^{\prime}(r)\right)^\prime=-\frac{1}{\lambda} r^{N-1}\tilde u_\alpha^p(r)\geq-\frac{1}{\lambda}r^{N-1}\,.
\end{equation}
Integrating between $\tilde\tau_\alpha$ and $r\in(\tilde\tau_\alpha,\tilde\sigma_\alpha)$, we get
$$
r^{N-1}\tilde u_\alpha^{\prime}(r) \geq-\frac{1}{\lambda N}(r^N-\tilde\tau_\alpha^N),
$$
which yields
\begin{equation}\label{eq5}
\tilde u_\alpha'(r)\geq-\frac{1}{\lambda N}r\qquad \text{for $r\in(\tilde\tau_\alpha,\tilde\sigma_\alpha)$}.
\end{equation}
Integrating again between $\tilde\tau_\alpha$ and $r\in(\tilde\tau_\alpha,\tilde\sigma_\alpha)$ and taking into account that $\tilde u_\alpha(\tilde\tau_\alpha)=1$, we have
\begin{equation}\label{eq9}
\tilde u_\alpha(r)\geq1-\frac{1}{2\lambda N}r^2\qquad \text{for $r\in(\tilde\tau_\alpha,\tilde\sigma_\alpha)$}.
\end{equation}
Therefore, for any fixed $r\in (0,l_1)$, letting $\alpha \to 0^+$ we infer that 
$$
1\geq \tilde u(r)\geq1-\frac{1}{2\lambda N}r^2,
$$
which implies that
\begin{equation}\label{eq10}
\lim_{r\to0}\tilde u(r)=1.
\end{equation}
Hence $\tilde u$ can be extended by continuity in $0$ by the value $\tilde u(0)=1$.  %As before integrating the equation 
%$$
%-r^{N-1} \tilde u_\alpha^\prime(r)%=\frac{1}{\lambda}\int_{\sigma_\e}^r\tilde v_\e^{p_\e}t^{N-1}\,dt\leq\frac{1}{\lambda}\int_{\sigma_\e}^rt^{N-1}\,dt.
%\leq\frac{1}{\lambda N}r^N,
%$$  
Moreover,  since $ \tilde u_\alpha^\prime(r)\leq0$ in $(\tilde\tau_\alpha, \tilde \sigma_\alpha)$, from \eqref{eq5} we get that
$$
| \tilde u_\alpha^\prime(r)|\leq\frac{r}{\lambda N}\qquad\text{for any $r\in(\tilde\tau_\alpha, \tilde \sigma_\alpha)$}.
$$
Therefore, fixing $r \in (0,l_1)$ and passing to the limit as $\alpha\to0^+$, we get
\begin{equation}\label{eq11}
|\tilde u'(r)|\leq\frac{r}{\lambda N}
\end{equation}
which gives $\lim_{r\to0}\tilde u'(r)=0$. 
Thus, $\tilde{u}$ in particular satisfies
$$
\left\{
\begin{array}{c}
 \tilde{u}''+\frac{N-1}{r}\tilde{u}'=-\frac{\tilde{u}^p}{\lambda} \quad  \hbox{ in } \ (0,l_1)\, ,\\[1ex]
 \tilde{u}(0)=1, \ \ \tilde{u}^\prime(0)=0\, , 
 \end{array}
 \right.
$$
and this implies that it  is of class $C^2$ up to 0. 

Summing up, $\tilde u$ is a smooth nontrivial radial solution to
$$
\begin{cases}
-\Mp(D^2 w)=w^{p} & \hbox{in} \ B_k,\\
\qquad  \ \ \ \ \ w> 0 & \hbox{in} \  B_k,\\
\qquad  \ \ \ \ \ w=0 & \hbox{on} \ \partial B_k.
\end{cases}
$$
Since $p>p_+^*$, this is a contradiction to the results of \cite{DLS, FQ}. The proof of Step 3 is complete.\\

\noindent\textbf{Conclusion}: By the previous steps, we deduce that, up to a  sequence  $\alpha\to 0^+$, $\tilde u_\alpha \to \tilde u$ in $C^2_{loc}((0,+\infty))$, where $\tilde u$ can be extended at the origin as a smooth positive radial solution of 
\begin{equation}\label{eq:probLimitRNMpm}
-\Mp (D^2 u)=u^p \quad\text{in $\R^N$.}
\end{equation}
Moreover, by construction, $\tilde u$ satisfies $\tilde{u}(0)=1$, it is radially decreasing in $[0,+\infty)$ and it is convex in $[l,+\infty)$. 

 If $p \in ]p^*_+ ,\frac{\Np+2}{\Np-2}]$, then it is known (see \cite{FQ}) that the only positive radial solutions of \eqref{eq:probLimitRNMpm} are pseudo-slow decaying solutions, which are functions changing concavity infinite times (see the proof of Lemma \ref{gammaalfa} (ix), or the proof of Proposition 3.3 in \cite{FQ}, and Remark  \ref{entire1}).  This clearly gives a contradiction because $\tilde u$ changes concavity exactly once.
We point out that this case is considered only for $\F=\Mp$, since for $\F=\Mm$ we have $p^*_->\frac{\Nm+2}{\Nm-2}$.

On the other hand, if $p>\frac{\Np+2}{\Np-2}$, again by the results of \cite{FQ}  the only positive radial solutions of \eqref{eq:probLimitRNMpm} are slow decaying functions,
  i.e. solutions for which the trajectory $\tilde{\gamma}(t)=(\tilde x(t),\tilde x'(t))$, associated with the Emden-Fowler transform $\tilde{x}(t)$ defined trough \eqref{EF}, satisfies
 \begin{equation}\label{asymgamma}
  \lim_{t\to +\infty} \tilde{\gamma} (t)= (c^*, 0)\, ,
\end{equation}
 where $c^*=(\Lambda\, \tilde{b}_+)^{\frac{1}{p-1}}$, with $\tilde{b}_+$ given by \eqref{coeff}.\\
Let us denote by $\tilde{x}_\alpha(t)$ the Emden-Fowler transform \eqref{EF} of $\tilde{u}_\alpha(r)$,  and let $\tilde{\gamma}_\alpha(t)=\left( \tilde{x}_\alpha (t), \tilde{x}_\alpha'(t)\right)$ be the associated trajectory, both   defined for $t\in \left[ \frac{p-1}{2} \log m_\alpha , \log \tilde{\rho}_\alpha\right]$. Then,  $\tilde{\gamma}_\alpha \to \tilde{\gamma}$ in $C^1_{loc}(\R)$.
We observe that for $t\in[\log \tilde\sigma_\alpha,\log \tilde\rho_\alpha ]$, the trajectory $\tilde{\gamma}_\alpha(t)$ belongs to $\mathcal{R}^-$ by Lemma \ref{gammaalfa} (v), and, therefore, the energy function $E$ defined in \eqref{energy} satisfies, by \eqref{Xeq},
$$
\frac{d}{dt} \left( E(\tilde{\gamma}_\alpha (t))\right) = \tilde{a}_+ \tilde{x}_\alpha'^2\, .
$$
Now, noticing that the condition $p>\frac{\Np+2}{\Np-2}$ is equivalent to  $\tilde{a}_+<0$, we infer that  $E(\tilde{\gamma}_\alpha (t))$ is monotone decreasing for $t \in [\log \tilde\sigma_\alpha, \log \tilde\rho_\alpha ]$. Hence,
$$
E (\tilde{\gamma}_\alpha (t)) \geq E (\tilde{\gamma}_\alpha (\log \tilde{\rho}_\alpha))= \frac{[\tilde x_\alpha^\prime(\log \tilde\rho_\alpha )]^2}{2} > 0\, .
$$
Letting $\alpha\to 0^+$, we deduce
$$
E(\tilde{\gamma}(t))\geq 0\quad \hbox{ for all } t\geq \log l\, .
$$
On the other hand, by \eqref{asymgamma}, we also have
$$
\lim_{t\to +\infty} E(\tilde{\gamma}(t))= E(c^*, 0)= - \frac{ (\Lambda\, \tilde{b}_+)^{\frac{2}{p-1}} \tilde{b}_+ (p-1)}{2\, (p+1)}<0\, ,
$$
and we reach a contradiction in this case as well.

\end{proof}

\section{Non-existence in the subcritical case $p<p^+_\pm$ for $\Mpm$}
In this and the next two sections we are concerned with the non-existence part of Theorem \ref{main}.

\begin{theorem}\label{Nonexisub}
Problem  \eqref{eq:probMpm} does not have  any positive radial  solution when $\F=\Mp$ and $1<p<p^*_+$, or when $\F=\Mm$ and $1<p<p^*_-$.
\end{theorem}
\begin{proof}
In view of \cite{AS}, it is sufficient to prove the result for $p>\frac{N_\pm}{N_\pm-2}$.
Let us consider $\F=\Mp$, the proof for $\F=\Mm$ being completely analogous. \\
We will show  that for any $\alpha>0$ the unique maximal solution $u_\alpha=u_\alpha(r)$ of \eqref{eq:initvalprob+} vanishes at some $\rho_\alpha\in (1,+\infty)$.
Consider the unique positive radial solution $u=u(r)$ of 
\begin{equation}\label{eq:probPallaUnit}
\begin{cases}
-\Mp(D^2 u)=u^{p} & \hbox{in} \ B,\\
%\qquad  \ \ \ \ \ u> 0 & \hbox{in} \  B\\
\qquad  \ \ \ \ \ u=0 & \hbox{on} \ \partial B,
\end{cases}
\end{equation}
whose existence is guaranteed by the assumption $p<p^*_+$  (see \cite[Theorem 5.1]{FQ}). 
Performing the Emden-Fowler change of variable \eqref{EF} for $u$, we obtain in the phase plane a  trajectory $\gamma(t)=(x(t),x'(t))$ defined for all $t\leq 0$ and satisfying 
$$
\lim_{t\to-\infty}\gamma(t)=(0,0)\, ,\quad \gamma(0)=(0,u'(1))\, .
$$
By Lemma \ref{gammaalfa} (v),   there exists a unique time $t_0<0$ such that $\gamma (t_0)=(x(t_0), 0)$ and $x(t_0)> c^*=(\Lambda\, \tilde{b}_+)^{\frac{1}{p-1}}$, $(c^*,0)$ being the unique equilibrium point.

We then deduce that the trajectory $\gamma(t)$ and the $x'$-axis bound a closed region containing the equilibrium points $(0,0)$ and $(c^*,0)$. 
Consider now the trajectory  $\gamma_\alpha(t)=(x_\alpha(t),x_\alpha'(t))$  relative to $u_\alpha$. Clearly, $\gamma_\alpha$ cannot intersect $\gamma$. Moreover it cannot approach neither the equilibrium points nor a periodic orbit, since otherwise it should cross $\gamma$. By Lemma \ref{gammaalfa} (vii), it follows that $\gamma_\alpha (t)$ must leave the fourth quadrant in finite time, i.e. there exists $\rho_\alpha\in(1,+\infty)$ such that $u_\alpha(\rho_\alpha)=0$ as desired.  
\end{proof}
%%%%%%%%%%%%%%%%%%%%%%%%%%%%%%%%%%%%%%%%%%%%%%%%%%%%%%%%%%%%%%%

\section{Non-existence at the critical level $p=p^*_+$ for $\Mp$}

\begin{theorem}\label{mainteo:NonExistenceMP}
If $\F=\Mp$ and $p=p^*_+$ then \eqref{eq:probMpm} does not have any positive radial  solution.
\end{theorem}
\begin{proof}
Let us assume  by contradiction that \eqref{eq:probMpm} admits a positive radial solution. Then there exists $u_\alpha=u_\alpha(r)$ solution of \eqref{eq:initvalprob+} defined in $[1,+\infty)$, and the trajectory $\gamma_\alpha(t):=(x_\alpha(t), x_\alpha^\prime(t))$,  associated with the Emden-Fowler transform \eqref{EF} of $u_\alpha$, is defined for $t\in [0,+\infty)$.\\
Let us also consider the trajectory $\gamma (t)=(x(t),x'(t))$, $t\in(-\infty,+\infty)$, associated with the unique (up to scaling) solution of
\begin{equation*}
-\Mp(D^2 u)=u^{p^*_+}  \quad\text{in $\R^N$}.
\end{equation*}
By \cite{FQ} (see also Remark \ref{entire1}), we know that  $\lim_{t\to\pm\infty}\gamma (t)=(0,0)$ and the support of $\gamma (t)$ bounds a compact region including the equilibrium point $(c^*,0)$ of the dynamical system \eqref{Xeq} for $p=p^*_+$, with $c^*=(\Lambda\, \tilde{b}_+)^{\frac{1}{p^*_+-1}}$. Since $\gamma_\alpha$ and $\gamma$ cannot intersect (apart from the origin),  we infer that $\gamma_\alpha(t)$ cannot  approach as $t\to +\infty$ neither a periodic orbit nor the equilibrium point $(c^*,0)$. Thus, from Lemma \ref{gammaalfa} (vii), it follows that $\gamma_\alpha(t)\to(0,0)$ as $t\to+\infty$ as well. Moreover,  Lemma \ref{gammaalfa} (v), which  also applies  to $\gamma$,   yields that both $\gamma_\alpha(t)$ and $\gamma(t)$  intersect (and cross) the $x$-axis exactly once. Hence, for $\varepsilon>0$ small enough, the vertical line $L_\e =\{ x=\varepsilon\}$  intersects exactly once $\gamma_\alpha$ and $\gamma$ in the fourth quadrant.  Therefore,   there exist unique $t_\varepsilon$ and $t_{\alpha,\e}$ such that 
$$
x(t_\varepsilon)=\e=x_\alpha(t_{\alpha,\e})\, , \quad 
x'_\alpha(t_{\alpha, \e})< x'(t_\e) <0 \, .
$$ 
Moreover, again by Lemma \ref{gammaalfa} (v), we can assume that $\varepsilon>0$ is so small that $\gamma(t)\,,\ \gamma_\alpha(t)\in \mathcal{R}^-$ respectively for $t\geq t_\e$ and $t\geq t_{\alpha, \e}$, see the picture below.
\begin{center}
\includegraphics[scale=1.3]{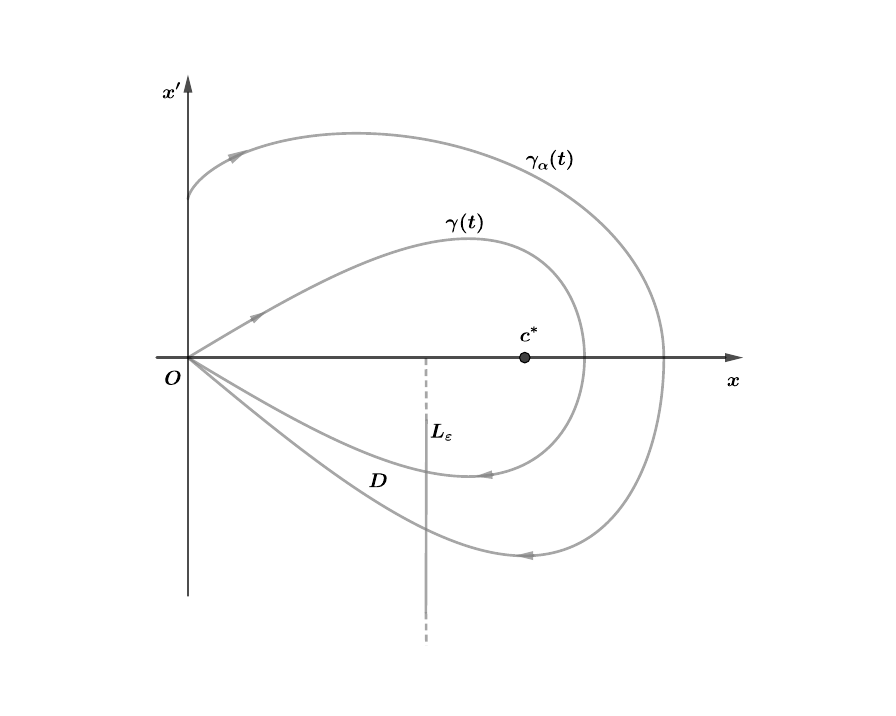}
\end{center}
Consider now  the compact region $D$ bounded by the closed  piecewise smooth curve formed by the trajectories $\gamma(t)$ for $t\geq t_\varepsilon$, $\gamma_{\alpha}(t)$ for $t\geq t_{\alpha,\e}$ and the vertical segment joining $\gamma(t_\varepsilon)$ and $\gamma_\alpha(t_{\alpha,\e})$.
Denoting by $|D|$ the area of $D$ and applying the Gauss-Green Theorem (with $\partial D$ clock-wise oriented),  we infer that
\beq\label{eq1teo2}
|D| = \int_{\partial D^{-}} y \ dx = \int_{t_{\alpha,\e}}^{+\infty} [x'_\alpha(t)]^2 \ dt + \int_{+\infty}^{t_\varepsilon} [x'(t)]^2 \ dt.
\eeq
On the other hand, considering the energy function $E$ defined in \eqref{energy}, we obtain that
$$  [x_\alpha^\prime(t)]^2= \frac{1}{\tilde{a}_+} \frac{d}{dt}\left[ E(\gamma_\alpha(t))\right] \,,  \quad   [x^\prime(t)]^2=\frac{1}{\tilde{a}_+} \frac{d}{dt}\left[ E(\gamma(t))\right] \, .
$$
Hence, using also that $x(t_\varepsilon)=x_\alpha(t_{\alpha,\e})$, we deduce
\begin{equation}\label{eq2teo2}
\left|D\right|=\frac{1}{\tilde a_+} \left[ E(\gamma(t_\e) -E(\gamma_\alpha(t_{\alpha, \e})\right]= \frac{1}{2 \tilde a_+}
\left[ \left[x^\prime(t_\varepsilon)\right]^2 -\left[x_\alpha^\prime(t_{\alpha,\e})\right]^2   \right].
\end{equation}
Now, since $p_+^* < \frac{\Np+2}{\Np-2}$, it holds that $\tilde a_+ =\tilde{a}_+(p^*_+)>0$, and since   $x_\alpha^\prime(t_{\alpha,\e}) < x^\prime(t_\varepsilon) < 0$,  from \eqref{eq2teo2} if  then follows that $|D| < 0$, which  clearly is a contradiction. 

\end{proof}

\begin{remark}
{\rm The previous proof does not work for $\F=\Mm$ and $p=p^*_-$ because in this case $p^*_- > \frac{\Nm+2}{\Nm-2}$ and thus $\tilde a _-(p^*_-)<0$. In particular \eqref{eq2teo2} does not give a contradiction.}
\end{remark}

\section{Non-existence at the critical level $p=p^*_-$ for $\Mm$}

\begin{theorem}\label{mainteo:NonExistenceMm}
If $\F=\Mm$ and $p=p^*_-$ then \eqref{eq:probMpm} does not have any positive radial solution.
\end{theorem}

The proof of Theorem \ref{mainteo:NonExistenceMm} will be given at the end of this section after some preliminary results on solutions of semilinear ODEs.\\

Let $\nu\in(2,+\infty)$. For any $p>1$ and $\alpha>0$, let $u_\alpha=u_\alpha(r)$ be the positive maximal solution of the initial value problem
\begin{equation}\label{4eq1}
\left\{\begin{array}{cl}
u''(r)+(\nu-1)\frac{u'(r)}{r}=- \frac{u^p(r)}{\lambda}  & \text{for $r>1$}\\[1ex]
u(1)=0,\;\,u'(1)=\alpha\,. &
\end{array}\right.
\end{equation}
The function $u_\alpha$ is defined on a maximal interval $[1,\rho_\alpha)$ with $\rho_\alpha\leq+\infty$. If $\rho_\alpha<+\infty$ then $u_\alpha(\rho_\alpha)=0$, otherwise $u_\alpha(r)>0$ for any $r>1$ and $u_\alpha(r)\to0$ as $r\to+\infty$.\\
Set $$D=\left\{\alpha\in(0,+\infty)\,:\;\rho_\alpha<+\infty\right\}.$$ The set $D$ is nonempty (see e.g. \cite[Proposition 3.2]{GLP}) and open by the continuous dependence on the initial data for \eqref{4eq1}. Let
\begin{equation}\label{4eq2}
\alpha^*=\alpha^*(p):=\inf D.
\end{equation}

In the next proposition, the behavior as $r\to +\infty$ of the solutions $u_\alpha$ is analyzed. In particular, when $\nu=N$ is an integer,  we recover the well known results (see e.g. \cite{JPY}) about  radial solutions of the semilinear exterior Dirichlet problem 
$$
\begin{cases}
-\Delta u=u^{p} & \hbox{in} \ \R^N\setminus \overline B,\\
  \ \ \ \ \ u=0 & \hbox{on} \ \partial B\, ,
\end{cases}
$$
in the case they exist, i.e.  in the supercritical regime $p>\frac{N+2}{N-2}$. For the sake of completeness, we include a (different) proof. Let us preliminarily observe that Lemma \ref{gammaalfa} clearly applies to the trajectory $\gamma_\alpha(t)=(x_\alpha(t), x'_\alpha(t))$, where $x_\alpha(t)$ is the Emden-Fowler transform \eqref{EF} of $u_\alpha(r)$, since problem \eqref{4eq1} is a special case of \eqref{eq:initvalprob-} with $\Lambda=\lambda$ and $N=\nu$.

\begin{proposition}\label{4prop1}
Let $p>\frac{\nu+2}{\nu-2}$. Then
\begin{itemize}
	\item[(i)] $\alpha^*>0$;\\
	\item[(ii)] if $\alpha>\alpha^*$ then $\rho_\alpha<+\infty$;\\
	\item[(iii)] $\rho_{\alpha^*}=+\infty$ and $\displaystyle\lim_{r\to+\infty}r^{\nu-2}u_{\alpha^*}(r)=C$, for some positive constant $C$;\\
	\item[(iv)] if $\alpha<\alpha^*$ then $\rho_{\alpha}=+\infty$ and $\displaystyle\lim_{r\to+\infty}r^{\frac{2}{p-1}}u_{\alpha}(r)=\left(\frac{2 \lambda ((\nu-2)p-\nu)}{(p-1)^2}\right)^{\frac{1}{p-1}}$.
\end{itemize}
\end{proposition}
\begin{proof}
(i). By contradiction let us assume $\alpha^*=0$. Then there is a sequence $\alpha_n\in D$ such that $\alpha_n\to0$ and $\rho_{\alpha_n}\to + \infty$ as $n\to+\infty$, in view of \cite[Lemma 3.1]{GLP}. \\
Let us  consider the   rescaled function
$$
\tilde u_{\alpha_n}(r)=\frac{1}{m_{\alpha_n}}u_{\alpha_n}\left(\frac{r}{m^{\frac{p-1}{2}}_{\alpha_n}}\right),\;\quad r>m^{\frac{p-1}{2}}_{\alpha_n},
$$
with $m_{\alpha_n}=u_{\alpha_n}(\tau_{\alpha_n})=\max_{[1,\rho_{\alpha_n}]}u_{\alpha_n}$.  Following the arguments of Step 1 and Step 2 of Theorem \ref{mainteo:Existence} we infer that $\tilde\tau_{\alpha_n}=m^{\frac{p-1}{2}}_{\alpha_n}\tau_{\alpha_n}\to0$   and $\tilde\rho_{\alpha_n}=m^{\frac{p-1}{2}}_{\alpha_n}\rho_{\alpha_n}$ is bounded from below by a positive constant.
Let us show that $\tilde\rho_{\alpha_n}\to+\infty$. If not, up to  a subsequence, $\tilde\rho_{\alpha_n}\to k\in(0,+\infty)$ and, again as in the  proof of Theorem \ref{mainteo:Existence}, $\tilde u_{\alpha_n}\to\tilde u$ in $C^2_{loc}(0,k)$ with 
\begin{equation}\label{eq:ODEpohoz}
\left\{\begin{array}{rl}
\tilde u''(r)+(\nu-1)\frac{\tilde u'(r)}{r}=-\frac{\tilde u^p(r)}{\lambda}  & \text{for $r>0$}\\[1ex]
\tilde u(0)=1,\;\,\tilde u'(0)=0,\;\, \tilde u(k)=0\, .
\end{array}\right.
\end{equation}
Multiplying the ODE in \eqref{eq:ODEpohoz} by $\tilde u(r)r^{\nu-1}$ we have 
$$(\tilde u^\prime(r)r^{\nu-1})^\prime \tilde u(r)+\frac{1}{\lambda}\tilde u^{p+1}(r)r^{\nu-1}=0\, .$$
Then, integrating by parts from  $0$ to $r$, this yields
\begin{equation}\label{eq:proPoho1}
\int_0^r (\tilde u^\prime)^2 s^{\nu-1} \ ds = \frac{1}{\lambda} \int_0^r \tilde u^{p+1} s^{\nu-1} \ ds + \tilde u(r)\tilde u^\prime(r)r^{\nu-1}.
\end{equation}
Similarly, multiplying \eqref{eq:ODEpohoz} by $\tilde u^\prime(r) r^{\nu}$ and integrating, we obtain
\begin{equation}\label{eq:proPoho2}
\left(\frac{\nu-2}{2}\right)\int_0^r (\tilde u^\prime)^2 s^{\nu-1} \ ds =- \frac{(\tilde u^\prime(r))^2}{2}  r^{\nu} - \frac{ \tilde u^{p+1}(r)}{\lambda(p+1)} r^{\nu} + \frac{\nu}{\lambda(p+1)} \int_0^r \tilde u^{p+1} s^{\nu-1} \ ds. 
\end{equation}
From \eqref{eq:proPoho1} and \eqref{eq:proPoho2} we deduce the following Pohozaev-like identity 
\begin{equation*}
r^\nu\frac{(\tilde u^\prime(r))^2}{2}+r^\nu\frac{\tilde u^{p+1}(r)}{\lambda(p+1)}+\frac{\nu-2}{2}r^{\nu-1}\tilde u(r)\tilde u'(r)=
\frac1\lambda\left(\frac{\nu}{p+1}-\frac{\nu-2}{2}\right)\int_0^r s^{\nu-1}\tilde u^{p+1}\,ds,
\end{equation*}
which leads to a contradiction for $r=k$, since  $\frac{\nu}{p+1}-\frac{\nu-2}{2}<0$ in view of the assumption $p>\frac{\nu+2}{\nu-2}$.\\
Hence, $\tilde\rho_{\alpha_n}\to+\infty$ and $\tilde u_{\alpha_n}$  converges  in $C^2_{loc}((0,+\infty))$ to $\tilde u$, with
\begin{equation}\label{4eq3}
\left\{\begin{array}{rl}
\tilde u''(r)+(\nu-1)\frac{\tilde u'(r)}{r}=- \frac{\tilde u^p(r)}{\lambda}  & \text{for $r>0$}\\[1ex]
\tilde u(0)=1,\;\,\tilde u'(0)=0 &. 
\end{array}\right.
\end{equation}
Denoting, respectively, by $\tilde x_{\alpha_n}$, $\tilde x$ the Emden-Fowler transform \eqref{EF} of $\tilde u_{\alpha_n}$, $\tilde u$,  and by $\tilde{\gamma}_{\alpha_n}=\left(\tilde x_{\alpha_n}, \tilde x_{\alpha_n}'\right)$, $\tilde{\gamma}=\left( \tilde x, \tilde x'\right)$ the associated trajectories, 
in the phase plane $\tilde \gamma_{\alpha_n}(t)\to\tilde \gamma(t)$ as $n\to+\infty$, the convergence being $C^1_{loc}(\R)$. We observe that in this case  $\tilde x_{\alpha_n}$ and $\tilde x$ are solution, respectively for $t\in\left(\frac{p-1}{2} \log
 m_{\alpha_n}, \log \tilde\rho_{\alpha_n} \right)$ and $t\in\R$,  of 
\begin{equation}\label{onlyone}
x''(t)= ax'(t)+b x(t)-\frac{x^p(t)}{\lambda} \, ,
\end{equation}
where now $a=\frac{\nu+2 -(\nu-2)p}{p-1}$ and $b=\frac{2\left( (\nu-2)p-\nu\right)}{(p-1)^2}$. Since $p>\frac{\nu+2}{\nu-2}$, then $a<0$ and, by a direct computation,  the  energy function
\begin{equation}\label{energyone}
E(x,y)=\frac{y^2}{2}- b \frac{x^2}{2}+\frac{x^{p+1}}{\lambda\, (p+1)}
\end{equation}
decreases along $\tilde{\gamma}_{\alpha_n}$ and $\tilde{\gamma}$. Hence
$$
 E\left( \tilde{\gamma}_{\alpha_n}(t)\right) >   E\left( \tilde{\gamma}_{\alpha_n}(\log \tilde\rho_{\alpha_n})\right)=\frac{(\tilde{x}'_{\alpha_n})^2\left( \log \tilde\rho_{\alpha_n} \right)}{2} >0\qquad\forall\, t\in\left(\frac{p-1}{2} \log
 m_{\alpha_n}, \log \tilde\rho_{\alpha_n} \right)
$$
and 
$$
 E\left(\tilde \gamma (t)\right) <\lim_{t\to-\infty}  E\left(\tilde \gamma (t)\right)= E(0,0)= 0 \qquad\forall\ t\in(-\infty,+\infty).
$$
Since $E\left( \tilde{\gamma}_{\alpha_n}(t)\right) \to E\left(\tilde \gamma (t)\right)$ for any fixed $t\in\R$, this leads to a contradiction.

\medskip
\noindent
(ii). Let $\alpha>\alpha^*$. Take $\bar\alpha\in D$ such that $\bar\alpha\in(\alpha^*,\alpha)$ and consider in the phase plane the trajectory $\gamma_{\bar\alpha}(t)=(x_{\bar\alpha}(t),x'_{\bar\alpha}(t))$, where $x_{\bar\alpha}$ is the Emden-Fowler transform \eqref{EF} associated with $u_{\bar\alpha}$.   By construction, $\gamma_{\bar\alpha}$ is defined for $t\in[0,\log \rho_{\bar\alpha}]$ and satisfies  $\gamma_{\bar\alpha}(0)=(0,\bar\alpha)$,  $\gamma_{\bar\alpha}(\log \rho_{\bar\alpha})=(0,x'_{\bar\alpha}(\log \rho_{\bar\alpha} ))$, with $x'_{\bar\alpha}(\log \rho_{\bar\alpha})<0$, and $x_{\bar\alpha}(t)>0$ for $t\in(0,\log \rho_{\bar\alpha})$. By Lemma \ref{gammaalfa} (v),  there exists a unique $T\in (0,\log \rho_{\bar\alpha})$ such that $x'_{\bar\alpha}(T)=0$ and, moreover,  $x_{\bar\alpha}(T)>c^*$, where $c^*=(\lambda b)^{\frac{1}{p-1}}$ and $(c^*,0)$ is the only equilibrium point in the right half plane for equation \eqref{onlyone}.  

Let us consider now  the trajectory $\gamma_\alpha(t)=(x_\alpha(t),x'_\alpha(t))$. If, by contradiction, $\rho_\alpha=+\infty$, then $\gamma_\alpha$ satisfies Lemma \ref{gammaalfa} (vii). On the other hand, in this case equation \eqref{onlyone} does not admit periodic solutions, due to the monotonicity along trajectories of the energy function \eqref{energyone}. Hence, $\gamma_\alpha (t)$ must approach as $t\to \infty$ either $(0,0)$ or $(c^*, 0)$. In both cases, it intersects the trajectory $\gamma_{\bar \alpha}$, which is impossible since $\alpha > \bar \alpha$. Thus,  the only possibility is that  $\gamma_\alpha$ leaves the right half plane in finite time, i.e. $\rho_\alpha<+\infty$ as desired.

\medskip
\noindent
(iii). Since $\alpha^*>0$ by (i) and the set $D$ is open, then necessarily $\alpha^*\notin D$, namely $\rho_{\alpha^*}=+\infty$.

In view of (ii), for any $\alpha>\alpha^*$ the trajectory $\gamma_\alpha(t)=(x_\alpha(t),x'_\alpha(t))$ crosses the $x'$-axis  at $t=\log \rho_\alpha$, with $\log \rho_\alpha \to +\infty$ as $\alpha\to\alpha^*$ due to the locally uniform convergence of $\gamma_{\alpha}$ to $\gamma_{\alpha^*}$, $\gamma_{\alpha^*}$ being the trajectory related to $u_{\alpha^*}$. Moreover, by the monotonicity of the energy function \eqref{energyone} along trajectories, one  has $E\left( \gamma_\alpha(t)\right) >0$ for any $t\in[0,\log \rho_\alpha]$, so that  $E\left( \gamma_{\alpha^*}(t)\right) \geq0$ for any $t\geq0$. 

Arguing as in the proof of (ii), we have that, as $t\to +\infty$, the trajectory $\gamma_{\alpha^*}(t)$ converges either to $(0,0)$ or to $(c^*, 0)$, with $c^*= (\lambda b)^{\frac{1}{p-1}}$. Since $E(c^*,0)=-  \frac{(c^*)^2 b (p-1)}{2\, (p+1)}<0$, the latter case is excluded and, by Lemma \ref{gammaalfa} (viii),  the conclusion follows.

\medskip
\noindent
(iv). Let $\alpha<\alpha^*$. By \eqref{4eq2},  we have $\rho_\alpha=+\infty$. Then, arguing as above,   either $\gamma_\alpha(t)\to(0,0)$ or $\gamma_\alpha(t)\to(c^*,0)$ as $t\to +\infty$. 
By contradiction, let us suppose   that  $\gamma_\alpha(t)\to(0,0)$. By Lemma \ref{gammaalfa} (viii),  this means that  the maximal solution $u_\alpha$ is of \eqref{4eq1} satisfies $\lim_{r\to+\infty}r^{\nu-2}u_\alpha(r)=C$ and $\lim_{r\to \infty} r^{\nu-1} u'_\alpha(r)=-(\nu-2)C$,  for a positive constant $C$.\\
 Let us now consider  the Kelvin transform of $u_\alpha$, defined by
\begin{equation}\label{K}
\bar u_\alpha(r):=\begin{cases}
r^{2-\nu}u_\alpha\left(\frac{1}{r}\right) & \text{if $r\in(0,1]$}\\
C & \text{if $r=0$}.
\end{cases}
\end{equation}
By a straightforward computation, one checks that $\bar u_\alpha$ is a bounded solution of 
\begin{equation}\label{4eq6}
\left\{\begin{array}{c}
u''(r)+(\nu-1)\frac{u'(r)}{r}=- \frac{1}{\lambda}r^{p(\nu-2)-\nu-2}u^p(r) \quad  \text{for $r\in(0,1)$},\\[2pt]
u(1)=0\, .  
\end{array}\right.
\end{equation}
Moreover, one has $\frac{\bar u_\alpha'(r)}{r}\to0$ ar $r\to0^+$, so that $\bar u_\alpha$ can be extended to a solution of \eqref{4eq6} for $r\in[0,1)$. Indeed, using the identity
\begin{equation}\label{4eq7}
(r^{\nu-1}\bar u_\alpha')'=-\frac1\lambda r^{p(\nu-2)-3}\bar u_\alpha^{p}\;,\qquad r\in(0,1),
\end{equation}
we first infer that $r^{\nu-1}\bar u_\alpha'$ is monotone decreasing. Furthermore, by definition \eqref{K} and by the decaying property as $r\to +\infty$ of $u_\alpha(r), \ u'_\alpha(r)$, we have
\begin{equation*}
\lim_{r\to0^+}r^{\nu-1}\bar u_\alpha'(r)=-\lim_{r\to+\infty}\left(u_\alpha'(r)r+(\nu-2)u_\alpha(r)\right)=0\, .
\end{equation*}
 Integrating \eqref{4eq7}, then we get
$$
r^{\nu-1}\bar u_\alpha'(r)=-\frac1\lambda\int_0^rs^{p(\nu-2)-3}\bar u_\alpha^p(s)\,ds
$$
 and
\begin{equation*}
\left|\frac{\bar u_\alpha'(r)}{r}\right| = \frac{1}{\lambda r^\nu}\int_0^rs^{p(\nu-2)-3}\bar u_\alpha^p(s)\,ds
\leq\frac{\left\|\bar u_\alpha^p\right\|_\infty}{\lambda(p(\nu-2)-2)}r^{p(\nu-2)-\nu-2}\to0\quad\text{as $r\to0^+$}.
\end{equation*}
Then $\bar u_\alpha$ is solution of \eqref{4eq6} also for $r\to0^+$.\\
By exactly the same argument,  the Kelvin transform $\bar u_{\alpha^*}$ of $u_{\alpha^*}$   also is a  solution of  \eqref{4eq6} for $r\in[0,1)$. But the  proof of \cite[Proposition 5.2]{ET} shows  that \eqref{4eq6} has a unique  positive solution in $[0,1)$. Indeed,  the proof  of  \cite{ET} is purely based  on the  analysis of the ODE problem \eqref{4eq6},  and  the used arguments do not depend on the first order coefficient to be natural or a real number larger than 1. Hence, by  uniqueness, we have $u_\alpha(r)=u_{\alpha^*}(r)$ for all $r\geq1$ and this is a contradiction since $\alpha<\alpha^*$. 
\end{proof}

\medskip

\begin{proof}[Proof of Thorem \ref{mainteo:NonExistenceMm}]
Assume by contradiction  that  problem \eqref{eq:probMpm}, with $\F=\Mm$ and $p=p^*_-$, admits a radial solution $u_\alpha=u_\alpha(r)$ with $u'_\alpha(1)=\alpha>0$. As in the proof of Theorem \ref{mainteo:NonExistenceMP}, let us consider the unique (up to scaling) entire radial solution $u_*$ of $-\Mm(D^2u)=u^{p^*_-}$ in $\R^N$ (see \cite{FQ}). After the Emden-Fowler transform \eqref{EF}, in the phase plane $u_\alpha$ and $u_*$ correspond to the trajectories
$\gamma_\alpha(t)=(x_\alpha(t),x'_\alpha(t))\quad\text{and}\quad\gamma_*(t)=(x_*(t),x'_*(t)),$
defined respectively for $t\in[0,+\infty)$ and $t\in\R$. Note that $\gamma_\alpha(t)$ and $\gamma_*(t)$ cannot intersect and, as in the proof of Theorem \ref{mainteo:NonExistenceMP}, we obtain that both $\gamma_\alpha(t)$ and $\gamma_*(t)$ converge to $(0,0)$ as $t\to+\infty$.\\
By Lemma \ref{gammaalfa} (iii) and (v), we have that there exist a unique $s_\alpha>0$ and a unique $s_*\in \R$ such that  
$\gamma_\alpha(s_\alpha)\, ,\ \gamma_*(s_*)\in \mathcal{C}$. Moreover, $\gamma_\alpha(t)$ and $\gamma_*(t)$ cross transversally from above the curve $\mathcal{C}$   for $t=s_\alpha$ and $t=s_*$ respectively, with $x^{p^*_--1}_\alpha(s_\alpha)\, ,\ x^{p^*_--1}_*(s_*)>\frac{\Lambda (N-1)}{p^*_-}$, and,  by \eqref{Xeq} with $p=p^*_-$ (and $\lambda$ and $\Lambda$ interchanged), $x_\alpha(t)$ and $x_*(t)$ satisfy, respectively for $t\geq s_\alpha$ and $t\geq s_*$,
\begin{equation}\label{frozen}
x''(t)=\tilde a_- x'(t)+\tilde b_- x(t)- \frac{x^{p^*_-}(t)}{\lambda}\, ,
\end{equation}
where $\tilde a_-$,  $\tilde b_-$ are defined in \eqref{coeff}.
We recall that along any trajectory $(x(t),x'(t))$ associated with a solution $x$ of \eqref{frozen}, the energy function $E$ defined by 
\begin{equation}\label{energy-}
E(x,y)=\frac{y^2}{2}-\tilde{b}_-\frac{x^2}{2}+\frac{x^{p^*_-+1}}{\lambda\, (p^*_-+1)}
\end{equation}
is monotone decreasing, since $\tilde a_-(p^*_-)<0$. Thus, equation \eqref{frozen} does not admit periodic orbits and its  only  equilibrium points in the phase plane are $(0,0)$ and $(c^*,0)$, with $c^*=(\tilde b_- \lambda)^{\frac{1}{p^*_--1}}$.

 Let us now consider the maximal positive solution $\bar x_\alpha=\bar x_\alpha(t)$ of the semilinear initial value problem \begin{equation*}
\left\{
\begin{array}{c}
x''(t)=\tilde a_- x'(t)+\tilde b_- x(t)- \frac{x^{p^*_-}(t)}{\lambda}\, ,\\[1ex]
x(s_\alpha)=x_\alpha(s_\alpha)\, ,\quad 
 x'(s_\alpha)=x_\alpha'(s_\alpha).
\end{array}
\right.
\end{equation*}
Similarly, let $\bar x_*=\bar x_*(t)$ be the maximal positive solution of
\begin{equation*}
\left\{
\begin{array}{r}
x''(t)=\tilde a_- x'(t)+\tilde b_- x(t)- \frac{x^{p^*_-}(t)}{\lambda}\, ,\\[1ex]
x(s_*)=x_*(s_*)\, ,\quad 
x'(s_*)=x_*'(s_*).
\end{array}
\right.
\end{equation*}
In other words, we  look at the problems with initial conditions given by $\gamma_\alpha(s_\alpha)$ and $\gamma_*(s_*)$, and with an equation having fixed coefficients   given by  $\tilde a_-$, $\tilde b_-$ and $1/\lambda$. \\
Clearly, 
  $\bar\gamma_\alpha(t):=(\bar x_\alpha(t),\bar x'_\alpha(t))=\gamma_\alpha(t)$ for  $t\geq s_\alpha$ and $\bar\gamma_*(t):=(\bar x_*(t), \bar x'_*(t))=\gamma_*(t)$ for $t\geq s_*$. Let us analyze the behavior of $\bar\gamma_\alpha(t)$ and $\bar\gamma_*(t)$ for $t<s_\alpha$ and $t<s_*$ respectively. We will  show  that there exist finite times $\bar s_\alpha<s_\alpha$ and $\bar s_*<s_*$ such that $\bar \gamma_\alpha(t)$ and $\bar \gamma_*(t)$ reach the $x'$-axis for $t=\bar s_\alpha$ and $t=\bar s_*$ respectively.\\
  Let us consider only the trajectory $\bar \gamma_*(t)$, the arguments for $\bar \gamma_\alpha$ being exactly the same. \\
We observe that, by the transversality of the intersection with $\mathcal{C}$ at the point $\bar \gamma_*(s_*)$, the trajectory $\bar \gamma_*(t)$ belongs to the region $\mathcal{R}^+$ for  $t$   in a left neighborhood of $s_*$. Let us first show  that $\bar \gamma_*(t)\in \mathcal{R}^+$ for  $t<s_*$. Indeed, assume, by contradiction, that there exists $t_*<s_*$ such that $\bar \gamma_* (t_*)\in \mathcal{C}$ and $\bar \gamma_*(t)\in \mathcal{R}^+$ for $t\in (t_*,s_*)$. By Lemma \ref{gammaalfa} (iii) and (iv) applied to $\bar \gamma_*$, we obtain that  $\bar x^{p^*_--1}_*(t_*) \leq \frac{\lambda (N-1)}{p^*_-}$ and that either $\bar \gamma_*(t)\in \mathcal{R}^-$ for all $t<t_*$ or $\bar \gamma_*$ intersects again $\mathcal{C}$ at some point on the right of $\frac{\lambda (N-1)}{p^*_-}$. In any case, $\bar \gamma_*(t)$ is bounded as $t\to -\infty$,  and  the Poincar\'e-Bendixson Theorem 
 implies that either $\bar \gamma_*(t)\to (0,0)$ or $\bar  \gamma_*(t)\to (c^*,0)$ as $t\to -\infty$. In both cases we reach a contradiction, since at the equilibrium points the energy \eqref{energy-} is non positive, whereas   along $\bar \gamma_*(t)$ it is decreasing and it converges to $0$ as $t\to +\infty$. This proves the claim. \\
  In order to show that $\bar \gamma_*(t)$ leaves the right half-plane in a finite time, let us use relation \eqref{EF} to switch from $\bar x_*(t)$ to $\bar u_*(r)$, which is the maximal positive solution of the initial value problem
$$
\left\{
\begin{array}{c}
u''+(\tilde N_- -1)\frac{u'}{r}=-\frac{u^{p^*_-}}{\lambda}\\[1ex]
u(\sigma_*)=u_*(\sigma_*)\, ,\quad u'(\sigma_*)=u'_*(\sigma_*)
\end{array}
\right.
$$
where $\sigma_*=e^{s_*}$ is the only zero of $u_*''(r)$. The solution $\bar u_*$ is defined on a maximal interval $(\bar r_*, +\infty)$, with $0\leq \bar r_*<\sigma_*$.\\
The fact that $\bar \gamma_*(t)\in \mathcal{R}^+$ for all $t<s_*$ means that $\bar u_*$ is a concave function in $(\bar r_*, \sigma_*]$. Then, $\bar u_*$ is bounded. If $(\bar u_*)'(r)<0$ for all $r\in (\bar r_*, \sigma_*]$, then $(\bar u_*)'$ also is bounded, since it is decreasing, and, therefore, $\bar \gamma_*(t)$ is bounded for  $t<\bar s_*$. Moreover, for $t<\bar s_*$, $\bar \gamma_*(t)$  lays in the portion of $\mathcal{R}^+$ below the line $L$ and, again by Poincar\'e-Bendixson Theorem, it follows that $\bar \gamma_* (t)\to (0,0)$ as $t\to -\infty$, which is impossible. Therefore, there exists $\bar \tau_*\in (\bar r_*, \sigma_*)$ such that $(\bar u_*)'(\bar \tau_*)=0$ and,  since $(\bar u_*)''(\bar \tau_*)=- \bar u^{p^*_-}_*(\bar \tau_*)/\lambda <0$, it follows that $(\bar u_*)'(r)>0$ for $r\in (\bar r_*, \bar \tau_*)$. 
\\ By writing, as usual, the equation satisfied by  $\bar u_*$ in the form
\begin{equation}\label{div}
\left( r^{\tilde N_- -1} \bar u_*'\right)'(r)= - \frac{  r^{\tilde N_- -1} \bar u_*^{p^*_-}(r)}{\lambda} < 0\, ,
\end{equation}
we further notice  that if, $\bar r_*=0$, then, for any fixed $r_0\in (0,\bar \tau_*)$, one has
$$
(\bar u_*)'(r)\geq \frac{r_0^{\tilde N_- -1} (\bar u_*)'(r_0)}{r^{\tilde N_- -1}}\quad \hbox{ for all } r\in (0,r_0]\, ,
$$
and this contradicts the boundedness of $\bar u_*$. Thus $\bar r_*>0$, and from \eqref{div}, it necessarily follows that
$$
\bar u_*(\bar r_*)=0\, ,\quad (\bar u_*)' (\bar r_*)=\, : \bar \alpha_*\in (0,+\infty)\, .
$$
The same arguments show that  $\bar x_\alpha (t)$ also is the Emden-Fowler transform of a function $\bar u_\alpha(r)$ which is the maximal positive solution of
$$
\left\{
\begin{array}{c}
u''+(\tilde N_- -1)\frac{u'}{r}=-\frac{u^{p^*_-}}{\lambda}\quad \hbox{ for } r>\bar r\\[1ex]
u(\bar r)=0\, ,\quad u'(\bar r)= \bar \alpha 
\end{array}
\right.
$$
for some $\bar r, \bar \alpha >0$. Then, the two functions $\bar u_*$ and $\bar u_\alpha$, suitably rescaled, are two distinct fast decaying solutions of problem  \eqref{4eq1} with $\nu=\tilde N_-$ and $p=p^*_->\frac{\tilde N_-+2}{\tilde N_--2}$, in contradiction with Proposition \ref{4prop1}.
\end{proof}

%%%%%%%%%%%%%%%%%%%%%%%%%%%%%%%%%%%%%%%%%%%%%%%%%%%%%%%%%%%%%%%
\section{Existence and uniqueness of fast decaying solutions}

In this section we complete the proof of Theorem \ref{main}, by showing in particular the existence and the uniqueness of fast decaying solutions of problems \eqref{eq:probMpm} for $p>p^*_\pm$.
\smallskip

For any $\alpha>0$, let $u_\alpha$ be the unique maximal solution of  either  \eqref{eq:initvalprob+} or \eqref{eq:initvalprob-}, defined on the maximal interval $[1,\rho_\alpha)$ with $\rho_\alpha\leq+\infty$.  As for the proof of Proposition \ref{4prop1}, let us set
 $$D=\left\{\alpha\in(0,+\infty)\,:\;\rho_\alpha<+\infty\right\}\, ,$$
 and 
 $$ \alpha^*=\alpha^*(p)=\inf D\, .
 $$
 Arguing as in the proof of Proposition \ref{4prop1}, one has $D=(\alpha^*,+\infty)$. Moreover, 
   by Theorems \ref{mainteo:Existence}, \ref{Nonexisub}, \ref{mainteo:NonExistenceMP} and \ref{mainteo:NonExistenceMm},    we have that $\alpha^*(p)>0$ if and only if  $p>p^*_\pm$. 
   
Thus, for $p>p^*_\pm$ and $0<\alpha\leq \alpha^*$, one has $\rho_\alpha=+\infty$. Let $\gamma_\alpha(t)=(x_\alpha(t),x'_\alpha(t))$ be the trajectory associated with the Emden-Fowler transform $x_\alpha (t)$ of $u_\alpha$, given by \eqref{EF}.\\
By continuous dependence on the initial data,  $\gamma_\alpha(t)\to \gamma_{\alpha^*}(t)$ in $C^1_{loc}([0,+\infty))$ as $\alpha \to \alpha^*$, and, since $\gamma_\alpha(t)$ reaches the $x'$-axis  in a finite time  for $\alpha>\alpha^*$, it follows necessarily that $\gamma_{\alpha^*}(t)\to (0,0)$ as $t\to +\infty$. Indeed, by Lemma \ref{gammaalfa} and the locally uniform convergence, it follows that $s_\alpha\to s_{\alpha^*}$ and $\gamma_{\alpha^*}(t)\in \mathcal{R}^-$ for all $t>s_{\alpha^*}$. This excludes the convergence of $\gamma_{\alpha^*}(t)$ to a periodic orbit as $t\to +\infty$, and the convergence to the equilibrium point $(c^*,0)$ (here $c^*= (\tilde{b}_+\Lambda)^{1/(p-1)}$ if $\mathcal{F}=\Mp$ and $c^*= (\tilde{b}_-\lambda)^{1/(p-1)}$ if $\mathcal{F}=\Mm$) is excluded as well,  by  arguing as in the conclusion of the proof of Theorem \ref{mainteo:Existence} and taking into account also that  $\gamma_{\alpha^*}$ cannot intersect the trajectory $\gamma$ associated with entire solutions. Therefore, $u_{\alpha^*}(r)$ is a fast decaying solution.

Next, in order to prove that $u_{\alpha^*}$ is the only fast decaying solution, let us consider separately  the cases of $\Mp$ and $\Mm$. 

Assume first that $u_\alpha$ solves \eqref{eq:initvalprob+} and, by contradiction, suppose that $\gamma_\alpha(t)\to (0,0)$ as $t\to +\infty$ for some $\alpha<\alpha^*$. If $p^*_+<p<\frac{\tilde{N}_++2}{\tilde{N}_+-2}$, then we can apply the argument of the proof of Theorem \ref{mainteo:NonExistenceMP} to the two trajectories $\gamma_{\alpha^*}(t)$ and $\gamma_\alpha(t)$, reaching a contradiction. On the other hand, if $p=\frac{\tilde{N}_++2}{\tilde{N}_+-2}$, then the trajectory associated with any fast decaying solution lies for sufficiently large $t$ on the zero level set of the energy function $E$ given by \eqref{energy}.
Hence, $\gamma_{\alpha^*}(t)$ and $\gamma_\alpha(t)$ definitively coincide, again a contradiction. Finally, if $p>\frac{\tilde{N}_++2}{\tilde{N}_+-2}$, then we can argue as in the proof of Theorem \ref{mainteo:NonExistenceMm} and, thanks to Proposition \ref{4prop1} for $\nu=\tilde N_+$, we obtain a contradiction as well. 

For the operator $\Mm$, the uniqueness of the fast decaying solution follows by the same proof of  Theorem \ref{mainteo:NonExistenceMm}.

Therefore, in both cases, we obtain that, for $\alpha<\alpha^*$, the trajectories $\gamma_\alpha(t)$ cannot approach the origin as $t\to +\infty$. 
Furthermore, since, for any $\alpha$, $\gamma_\alpha(t)$ cannot intersect the trajectory $\gamma(t)$ associated with any entire solution $u$ of $-\Mpm(D^2u)=u^p$, and since,  for $p>p^*_\pm$ and $\alpha<\alpha^*$, neither $u$ nor $u_\alpha$ is a fast decaying solution, it follows that   $\gamma(t)$ and $\gamma_\alpha(t)$ must have the same behavior as $t\to +\infty$. 

According to Theorem 1.1 of \cite{FQ}, we then obtain the following result.
 \begin{theorem}\label{fastMp} Let $p>p^*_+$ and let $u_\alpha$ denote  the maximal solution of \eqref{eq:initvalprob+}. Then:
 \begin{itemize}
 \item[(i)]  $u_{\alpha^*}$  is a fast decaying solution;
 \item[(ii)] for any $\alpha<\alpha^*$, $u_\alpha$ is a pseudo-slow decaying solution   if $p^*_+<p\leq \frac{\tilde{N}_++2}{\tilde{N}_+-2}$;
  \item[(iii)] for any $\alpha<\alpha^*$, $u_\alpha$ is a  slow decaying solution  if $p>\frac{\tilde{N}_++2}{\tilde{N}_+-2}$.
  \end{itemize}
\end{theorem}
 Analogously, by applying Theorem 1.2 of \cite{FQ}, we deduce the following theorem.
 \begin{theorem}\label{fastMm} Let $p>p^*_-$ and let $u_\alpha$ denote  the maximal solution of \eqref{eq:initvalprob-}. Then:
 \begin{itemize}
 \item[(i)]  $u_{\alpha^*}$  is a fast decaying solution;
 \item[(ii)] for any $\alpha<\alpha^*$, $u_\alpha$ is either a pseudo-slow or a slow decaying solution   if $p^*_-<p\leq \frac{N+2}{N-2}$;
  \item[(iii)] for any $\alpha<\alpha^*$, $u_\alpha$ is a  slow decaying solution  if $p> \frac{N+2}{N-2}$.
  \end{itemize}
\end{theorem}

\end{document}